\documentclass{article}[utf8, 12pt]%
\usepackage{amsmath}
\usepackage{amsfonts}
\usepackage{amssymb}
\usepackage{mathrsfs}
\usepackage{graphicx}
\usepackage{hyperref}
\usepackage{stmaryrd}
\usepackage{authblk}
\usepackage{enumerate}
\usepackage{bbm}
\usepackage{color}%
\usepackage{esint}
\usepackage{mathtools}
\usepackage[title,titletoc]{appendix}
\usepackage{xpatch,amsthm}
\makeatletter
\xpatchcmd{\@thm}{\fontseries\mddefault\upshape}{}{}{} 
\makeatother
\providecommand{\U}[1]{\protect\rule{.1in}{.1in}}
\newtheorem{theorem}{Theorem}

\newtheorem{definition}[theorem]{Definition}

\newtheorem{lemma}[theorem]{Lemma}

\newtheorem*{question*}{Question}
\newtheorem{proposition}[theorem]{Proposition}
\newtheorem{remark}[theorem]{Remark}

\begin{document}
	\author[a]{Banhirup Sengupta \thanks{banhirups@iisc.ac.in}}
	\author[b]{Swarnendu Sil\thanks{swarnendusil@iisc.ac.in}}
	\affil[a,b]{Department of Mathematics \\ Indian Institute of Science\\ Bangalore, India}
	\title{Morrey-Lorentz estimates for Hodge-type systems}
	\maketitle 
	\begin{abstract}
		We prove up to the boundary regularity estimates in Morrey-Lorentz spaces for weak solutions of the linear system of differential forms with regular anisotropic  coefficients 
		\begin{equation*}
			d^{\ast} \left( A  d\omega \right)  + B^{\intercal}d d^{\ast} \left( B\omega \right)  =  \lambda B\omega + f  \text{ in } \Omega,
		\end{equation*}
		with either $ \nu\wedge \omega$ and $\nu\wedge d^{\ast} \left( B\omega \right)$ or $\nu\lrcorner B\omega$ and 
		$\nu\lrcorner \left( A d\omega \right)$ prescribed on $\partial\Omega.$ We derive these estimates from the $L^{p}$ estimates obtained in \cite{Sil_linearregularity} in the spirit of Campanato's method. Unlike Lorentz spaces, Morrey spaces are neither interpolation spaces nor rearrangement invariant. So Morrey estimates can not be obtained directly from the $L^{p}$ estimates using interpolation. We instead adapt an idea of Lieberman \cite{Lieberman_morrey_from_Lp} to our setting to derive the estimates. Applications to Hodge decomposition in Morrey-Lorentz spaces, Gaffney type inequalities and estimates for related systems such as Hodge-Maxwell systems and `div-curl' systems are discussed. 
	\end{abstract}
\tableofcontents \bigskip

	\textbf{Keywords:} Morrey estimates, Morrey-Lorentz estimates, Boundary regularity, differential form, Campanato method, Hodge Laplacian, Maxwell system. \\
	
	\textbf{MSC codes:} 35J57, 35B65, 35Q60. 
	\section{Introduction}
	Perhaps the most important second order elliptic systems for differential forms are the Poisson problem for the Hodge Laplacian with prescribed `tangential part' or prescribed `normal part' on the boundary respectively, namely the systems,
	\begin{equation*}
		\left\lbrace \begin{gathered}
			d^{\ast} d\omega  + dd^{\ast} \omega  =  f  \text{ in } \Omega, \\
			\nu\wedge \omega = 0 \text{  on } \partial\Omega. \\
			\nu\wedge d^{\ast} \omega = 0 \text{ on } \partial\Omega. 
		\end{gathered} 
		\right. \qquad \text{ or } \qquad \left\lbrace \begin{gathered}
			d^{\ast} d\omega  + dd^{\ast} \omega   =  f  \text{ in } \Omega, \\
			\nu\lrcorner \omega = 0 \text{  on } \partial\Omega. \\
			\nu\lrcorner d\omega = 0 \text{ on } \partial\Omega. 
		\end{gathered}\right.
	\end{equation*} 
	Standard $L^{p}$ and Schauder estimates for these systems are established by Morrey \cite{MorreyHarmonic2} (see also \cite{Morrey1966}). These estimates lead to the Hodge decompositions and a large number of related results ( see \cite{SchwarzHodge} for much more on this ). The proof of Morrey is based on representation formulas for the components of the solution using the Green's function for scalar Laplacian. This makes crucial use of the fact that as far as the principal order terms are concerned, the whole system decouples and gets reduced to $\binom{n}{k}$ number of scalar Poisson problems with lower order terms, out of which some has Dirichlet boundary condition and the rest has Neumann boundary condition. This method however, does not work for the general `Hodge systems' 
	\begin{equation*}
		\left\lbrace \begin{gathered}
			d^{\ast} \left( Ad\omega \right) + B^{\intercal}dd^{\ast}\left( B \omega \right)  =  f  \text{ in } \Omega, \\
			\nu\wedge \omega = 0 \text{  on } \partial\Omega. \\
			\nu\wedge d^{\ast} \left( B \omega \right)  = 0 \text{ on } \partial\Omega. 
		\end{gathered} 
		\right.  \text{ or }  \left\lbrace \begin{gathered}
			d^{\ast} \left( Ad\omega \right) + B^{\intercal}dd^{\ast}\left( B \omega \right)   =  f  \text{ in } \Omega, \\
			\nu\lrcorner \left( B \omega \right) = 0 \text{  on } \partial\Omega. \\
			\nu\lrcorner \left( Ad\omega \right) = 0 \text{ on } \partial\Omega. 
		\end{gathered}\right.
	\end{equation*}
	The presence of the matrices $A$, $B$ prevents a decoupling of the system. Verifying either Lopatinski\u{i}-Shapiro or 
	the Agmon-Douglis-Nirenberg conditions to show that these boundary value problems are elliptic is also prohibitively tedious. The systems, however, are important for applications, e.g. to time-harmonic Maxwells equations when the coefficient tensors are anisotropic.
	
	\par  Standard $L^{p}$ and Schauder estimates for these systems were established, using the Campanato-Stampacchia method, in \cite{Sil_linearregularity}. Since the method uses  interpolation to obtain the $L^{p}$ estimates, estimates in $L^{p}$-based Morrey spaces can not be obtained this way, as Morrey spaces are neither interpolation spaces nor rearrangement invariant. 
	
	\par In this article, we derive Morrey-Lorentz estimates for these systems, still in the spirit of Campanato's method, from the $L^{p}$ estimates in \cite{Sil_linearregularity}. The main idea is to prove and use suitable decay estimates for the Lorentz quasinorms. This technique goes back to Lieberman \cite{Lieberman_morrey_from_Lp}, who used this in the context of $L^{p}$ spaces to derive estimates in $L^{p}$-based Morrey spaces from a suitable form of the $L^{p}$ estimates for scalar elliptic equations with Dirichlet or Oblique derivative type boundary conditions. We modify the method to our setting, where specific features of the Hodge-type systems comes into play. On top of this, the generalization to Morrey-Lorentz spaces also adds some technical complication. This is due to the fact that though $L^{p}$ spaces are always separable and reflexive for $1 < p < \infty,$ the Lorentz spaces $\mathcal{L}^{\left(p, \theta\right)}$ are not reflexive for $\theta = 1 $ or $\infty$ and not separable for $\theta = \infty.$  This is reflected in our argument in Lemma \ref{approximation lemma}. 
	
	\par  As a consequence of our main estimates, we derive a host of results, namely Theorem \ref{generalHodgesystemtheorem}, Theorem \ref{generalHodgesystemtheoremnormal}, Theorem \ref{Hodge decomposition}, Theorem \ref{Maxwell in MorreyLorentz},  Theorem \ref{divcurl system MorreyLorentz}, Theorem \ref{gaffney}. To the best of our knowledge, these results are not only new in this generality, but also new for the pure Morrey case $p = \theta.$ All our results work for $\mathbb{R}^{N}$-valued differential forms as well. Our results extends the theory, started in \cite{Sil_linearregularity}, to the context of Morrey-Lorentz spaces. 
	\par The rest of our article is organized as follows. We record our notations and preliminaries about Morrey-Lorentz-Sobolev spaces in Section \ref{Preliminaries}. Section \ref{Sobolev inequalities in Morrey-Lorentz} proves the Sobolev and Poincar\'{e}-Sobolev type inequalities in these spaces which we would use.  Our main estimates are proved in Section \ref{crucial estimates}. Section \ref{main results} states and proves our main results. For the sake of clarity, we state and prove only second order estimates. However, our technique can easily be adapted to derive gradient estimates. 
	
	\section{Preliminaries}\label{Preliminaries}
	\subsection{Notations}\label{notations}
	We record the notations we would use for exterior forms. For further details we refer to
	\cite{CsatoDacKneuss} and \cite{silthesis}. Let $n \geq 2,$ $N \geq 1$ and $0 \leq k \leq n$ be integers. 
	\begin{itemize}
		\item  The vector space of all alternating $k-$linear maps
		$f:\underbrace{\mathbb{R}^{n}\times\cdots\times\mathbb{R}^{n}}_{k-\text{times}%
		}\rightarrow\mathbb{R}$ will be denoted by $\Lambda^{k}\mathbb{R}^{n},$ with $\Lambda^{0} \mathbb{R}%
		^{n}:=\mathbb{R}.$ For vector-valued forms, we introduce some shorthand. We denote 
		\begin{align*}
			\varLambda^{k}:= \Lambda^{k}\mathbb{R}^{n}\otimes \mathbb{R}^{N}.
		\end{align*}
		$\left\{  e_{1},\cdots,e_{n}\right\}$ is the standard basis of $\mathbb{R}%
		^{n}.$ The dual basis $\left\{  e^{1},\cdots,e^{n}\right\} $ is a basis for $\Lambda^{1}\mathbb{R}^{n}$ and $	\left\lbrace e^{I}:= e^{i_{1}}\wedge\cdots\wedge e^{i_{k}} \right\rbrace_{1\leq i_{1}<\cdots<i_{k}\leq n} $ 	is a basis of $\Lambda^{k}\mathbb{R}^{n}.$ An element $\xi\in\varLambda^{k}$ will therefore be written as
		\begin{align*}
			\xi=\sum\limits_{j=1}^{N}\sum_{I\in\mathcal{T}^{k}}\xi
			_{I,j}\,e^{I}\otimes e_{j}, \quad \mathcal{T}^{k}=\left\{ \left(  i_{1}\,,\cdots,i_{k}\right)
			:1\leq i_{1}<\cdots<i_{k}\leq n\right\}  .
		\end{align*} 
		
		\item $\wedge,$ $\lrcorner\,,$ $\left\langle \ ;\ \right\rangle $ and,
		respectively, $\ast$ denote the exterior product, the interior product, the
		scalar product and, respectively, the Hodge star operator, extended componentwise in the obvious fashion to vector-valued forms. 
	\end{itemize}
	Now we describe our notations for the sets we would be using a lot. 
	\begin{itemize}
		\item For any Lebesgue measurable subset $A \subset \mathbb{R}^{n}, $ we denote its $n$-dimensional Lebesgue measure  by $\left\lvert A \right\rvert.$ 
		\item For any $z \in \mathbb{R}^{n}$ and any $r>0,$ the open ball with center $z$ and radius $r$ is denoted by $B_{r}\left( z\right) := \left\lbrace x \in \mathbb{R}^{n}: \left\lvert x - z \right\rvert < r\right\rbrace$. We would just write $B_{r}$ when the center of the ball is the origin, i.e. when $z=0 \in \mathbb{R}^{n}.$
		\item The open upper half space is denoted by 
		\begin{align*}
			\mathbb{R}^{n}_{+} := \left\lbrace x = \left(x', x_{n}\right) \in \mathbb{R}^{n}: x_{n} >0\right\rbrace, 
		\end{align*} The boundary of the open upper half space is denoted as 
		\begin{align*}
			\partial\mathbb{R}^{n}_{+} := \left\lbrace x = \left(x', 0\right) \in \mathbb{R}^{n}: x' \in \mathbb{R}^{n-1}\right\rbrace. 
		\end{align*} For any $z \in \partial\mathbb{R}^{n}_{+}$ and any $r>0,$ the open upper half ball with center $z$ and radius $r$ are denoted by $B^{+}_{r}\left( z\right) := \left\lbrace x = \left(x', x_{n}\right) \in \mathbb{R}^{n}: \left\lvert x - z \right\rvert < r , x_{n} >0 \right\rbrace $.  
		We would just write $B^{+}_{r}$ when the center of the balls is the origin, i.e. when $z=0 \in \mathbb{R}^{n}.$ For us, $\Gamma_{r}\left( z\right)$  and $\Sigma_{r}\left( z\right)$ would denote the flat part and the curved part, respectively, of the boundary of the half ball $B_{r}^{+}\left( z\right).$ More precisely, 
		\begin{align*}
			\Gamma_{r}\left( z\right):= \partial B_{r}^{+}\left( z\right) \cap 	\partial\mathbb{R}^{n}_{+}\qquad  \text{ and } \qquad \Sigma_{r} := \partial B_{r}^{+}\left( z\right)\setminus \Gamma_{r}\left( z\right).  
		\end{align*}
		\item  For any  open subset $\Omega\subset\mathbb{R}^n,$ and for any $z \in \mathbb{R}^{n}$ and any $r>0,$ we denote 
		\begin{align*}
			\Omega_{\left(r, z\right)}:= B_{r}\left( z\right) \cap \Omega. 
		\end{align*}
		Once again, we would write $\Omega_{\left(r\right)}$ when $z=0 \in \mathbb{R}^{n}.$
		
		\item Let $\mathcal{U}\subset\mathbb{R}_+^n$ be a smooth open set which is star-shaped about the origin such that 
		\begin{align*}
			B^{+}_{1/2} \subset B^{+}_{3/4} \subset \mathcal{U} \subset B^{+}_{7/8} \subset B^{+}_{1}.  
		\end{align*}
		Note that this implies $\mathcal{U}$ is contractible and $	\Gamma_{3/4} \subset \partial\mathcal{U}.$	For any $x_{0} \in \partial \mathbb{R}^{n}_{+},$ we set 
		\begin{align*}
			\mathcal{U}_{R}\left(x_{0}\right) :=  \left\lbrace x_{0} + R x : x \in \mathcal{U} \right\rbrace = \left\lbrace x \in \mathbb{R}_{+}^{n}: \frac{1}{R}\left( x - x_{0}\right) \in \mathcal{U} \right\rbrace. 
		\end{align*}
		We also write $\mathcal{U}_{R} := \mathcal{U}_{R}\left(0\right).$
		
		\item For the rest, $\Omega\subset\mathbb{R}^n$ will always denote an open, bounded subset with at least Lipschitz boundary. $\nu$ will always denote the outward unit normal field to $\partial\Omega,$ which will be identified with the $1$-form 
		$\displaystyle \nu = \sum_{i=1}^{n} \nu_{i} dx^{i} .$
	\end{itemize}
	For any finite vector space $X$ over the reals, the notation $\operatorname{Hom}\left( X\right)$ would denote the vector space of linear maps $A:X\rightarrow X.$ $A^{\intercal}$ would denote the adjoint or transpose of $A \in \operatorname{Hom}\left( X\right).$ \begin{definition}\label{legendre condition}
		A bounded measurable map $A \in L^{\infty}\left( \Omega; \operatorname{Hom}(\varLambda^{k}) \right)$ is called \textbf{uniformly Legendre elliptic} if there exists a constant $ \gamma >0$ such that  we have 
		$$  \langle A (x)  \xi \  ;\  \xi \rangle  \geq  \gamma \left\vert  \xi \right\vert^{2}	\quad  \text{ for every } \xi \in \varLambda^{k} \text{ and for   a.e. }  x \in \Omega.$$ 
	\end{definition}
	Clearly, if $A \in L^{\infty}\left( \Omega; \operatorname{Hom}(\varLambda^{k}) \right)$ is uniformly Legendre elliptic, then $A^{\intercal} \in L^{\infty}\left( \Omega; \operatorname{Hom}(\varLambda^{k}) \right) $, defined as the matrix field $x \mapsto \left(A\left(x\right)\right)^{\intercal}$ is also uniformly Legendre elliptic. We would often just say $A$ satisfies the Legendre condition or is Legendre elliptic and $\gamma$ would always stand for the ellipticity constant. 
	\subsection{Function spaces for differential forms}	
	\begin{itemize}
		\item A $\mathbb{R}^{N}$-valued differential $k$-form $\omega$ on $\Omega$ is a
		measurable function $\omega:\Omega\rightarrow\varLambda^{k}.$ The usual Lebesgue, Sobolev and H\"{o}lder spaces and their local versions are defined componentwise in the usual way and are denoted by their usual symbols. Morrey-Lorentz spaces are defined in section \ref{function spaces}. 
		\item Two special differential operators on differential forms will have a special significance for us. A $\mathbb{R}^{N}$-valued differential
		$(k+1)$-form $\varphi\in L^{1}_{\text{loc}}(\Omega;\varLambda^{k+1})$
		is called the exterior derivative of $\omega\in
		L^{1}_{\text{loc}}\left(\Omega;\varLambda^{k}\right),$ denoted by $d\omega$,  if
		$$
		\int_{\Omega} \eta\wedge\varphi=(-1)^{n-k}\int_{\Omega} d\eta\wedge\omega,
		$$
		for all $\eta\in C^{\infty}_{c}\left(\Omega;\varLambda^{n-k-1}\right).$  The Hodge codifferential of $\omega\in L^{1}_{\text{loc}}\left(\Omega;\varLambda^{k}\right)$ is
		an $\mathbb{R}^{N}$-valued $(k-1)$-form, denoted $d^{\ast}\omega\in L^{1}_{\text{loc}}\left(\Omega;\varLambda^{k-1}\right)$
		defined as
		$$
		d^{\ast}\omega:=(-1)^{nk+1} \ast d \ast \omega.
		$$ See \cite{CsatoDacKneuss} and \cite{silthesis} for the properties and the integration by parts formula regarding these operators.
	\end{itemize}
	The spaces $W_{T}^{1,2}\left(  \Omega;\varLambda^{k}\right)  $ and $W_{N}^{1,2}\left(  \Omega;\varLambda^{k}\right)  $ are defined as ( see \cite{CsatoDacKneuss} )
	\begin{align*}
		W_{T}^{1,2}\left(  \Omega;\varLambda^{k}\right)  &=\left\{  \omega\in
		W^{1,2}\left(  \Omega;\varLambda^{k}\right)  :\nu\wedge\omega=0\text{ on
		}\partial\Omega\right\}, \\
		W_{N}^{1,2}\left(  \Omega;\varLambda^{k}\right)  &=\left\{  \omega\in
		W^{1,2}\left(  \Omega;\varLambda^{k}\right)  :\nu\lrcorner\omega=0\text{ on
		}\partial\Omega\right\} 
	\end{align*}
	The subspaces $W_{d^{\ast}, T}^{1,2}(\Omega; \varLambda^{k})$ and $W_{d, N}^{1,2}(\Omega; \varLambda^{k})$  are defined as  
	\begin{align*}
		W_{d^{\ast}, T}^{1,2}(\Omega; \varLambda^{k}) &= \left\lbrace \omega \in W_{T}^{1,2}(\Omega; \varLambda^{k}) : d^{\ast}\omega = 0 \text{ in }
		\Omega \right\rbrace, \\
		W_{d, N}^{1,2}(\Omega; \varLambda^{k}) &= \left\lbrace \omega \in W_{N}^{1,2}(\Omega; \varLambda^{k}) : d\omega = 0 \text{ in }
		\Omega \right\rbrace.		
	\end{align*}
	The space of tangential and normal harmonic $k$-fields are defined as 
	\begin{align*}
		\mathcal{H}^{k}_{T}\left(  \Omega;\varLambda^{k}\right)  &=\left\{  \omega\in
		W_{T}^{1,2}\left(  \Omega;\varLambda^{k}\right)  :d\omega=0\text{ and }%
		d^{\ast}\omega=0\text{ in }\Omega\right\},\\
		\mathcal{H}^{k}_{N}\left(  \Omega;\varLambda^{k}\right)  &=\left\{  \omega\in
		W_{N}^{1,2}\left(  \Omega;\varLambda^{k}\right)  :d\omega=0\text{ and }%
		d^{\ast}\omega=0\text{ in }\Omega\right\}
	\end{align*}
	For a given $B \in C^{l+2}\left(\overline{\Omega}; \operatorname{Hom}\left(\varLambda^{k}\right)\right)$, satisfying the Legendre condition, let us define the space 
	$$ W^{1,2}_{B,N}\left(\Omega; \Lambda^{k} \right):= \left\lbrace \omega \in W^{1,2}\left(\Omega; \Lambda^{k} \right): \nu\lrcorner \left( B(x)\omega \right) = 0 \text{ on } 
	\partial\Omega \right\rbrace . $$
	For half-balls, we need the following subspaces. 
	\begin{align*}
		W_{T, \text{flat}}^{1,2}&(B_{r}^{+}\left(z\right) ; \varLambda^{k}) \\&= \left\lbrace \psi \in W^{1,2}(B_{r}^{+}\left(z\right) ; \varLambda^{k}):
		e_n \wedge \psi = 0 \text{ on }  \Gamma_{r}\left(z\right), \psi = 0 \text{ near } \Sigma_{r}\left(z\right) \right\rbrace, \\ 
		W_{N, \text{flat}}^{1,2}&(B_{r}^{+}\left(z\right) ; \varLambda^{k}) \\&= \left\lbrace \psi \in W^{1,2}(B_{r}^{+}\left(z\right) ; \varLambda^{k}):
		e_n \lrcorner \psi = 0 \text{ on }  \Gamma_{r}\left(z\right), \psi = 0 \text{ near } \Sigma_{r}\left(z\right) \right\rbrace. 
	\end{align*}
	Here $\nu \wedge \omega$ and $\nu \lrcorner \omega$ denotes the tangential and normal trace, respectively, on the boundary. A crucial fact about these traces that we would constantly use is the following ( see \cite{CsatoDacKneuss}, \cite{SchwarzHodge}, \cite{Morrey1966} ). 
	\begin{proposition}\label{identity for traces}
		Let $u \in W^{1,p}\left(  \Omega;\varLambda^{k}\right)$ for any $1 < p < \infty.$ Then 
		\begin{align*}
			\nu \wedge \omega = 0  \quad \text{ on } \partial\Omega \qquad &\Rightarrow \qquad \nu \wedge d\omega = 0 \quad \text{ on } \partial\Omega, \\
			\nu \lrcorner \omega = 0  \quad \text{ on } \partial\Omega  \qquad &\Rightarrow \qquad \nu \lrcorner d^{\ast}\omega = 0 \quad \text{ on } \partial\Omega. 
		\end{align*} 
	\end{proposition}
	\subsection{Morrey and Lorentz type spaces}\label{function spaces}
	\subsubsection{Morrey-Lorentz spaces}
	\begin{definition}[Morrey Spaces]
		Let $1\leqslant p <  \infty $ and $ 0 \leq \mu < n$ be real numbers. The Morrey space $\mathrm{L}^{p}_{\mu}\left(\Omega\right) $ stands for the space of all $f \in L^{p}\left(\Omega\right)$ such that 
		$$ \lVert f \rVert_{\mathrm{L}^{p}_{\mu}\left(\Omega\right)}^{p} := \sup_{\substack{ x_{0} \in \overline{\Omega},\\ \rho >0 }} 
		\frac{1}{\rho^{\mu}} \int_{\Omega_{(\rho , x_{0})}} \lvert f \rvert^{p} < +\infty, $$ endowed with the norm 
		$ \lVert f \rVert_{\mathrm{L}^{p}_{\mu}\left(\Omega\right)}.$
	\end{definition} 
	Morrey spaces were introduced by Morrey in \cite{Morrey_MorreySpacesintro}. Now we define the Lorentz spaces, introduced by Lorentz in \cite{Lorentz_LorentzSpacesIntroFirst}. 
	\begin{definition}[Lorentz Spaces]
		Let $1\leqslant p <  \infty$ and $1 \leq \theta < \infty$ be real numbers. A measurable function $f:\Omega \rightarrow \mathbb{R}$ is said to belong to the Lorentz space $\mathcal{L}^{(p,\theta)} \left(\Omega\right)$ if 
		\begin{align*}
			\left\lVert f \right\rVert_{\mathcal{L}^{(p,\theta)}\left(\Omega\right)} 
			:= \left( \int_{0}^{\infty} \left( t^{\frac{1}{p}} f^{\ast}_{\Omega}\left(t\right)\right)^{\theta} \frac{\mathrm{d}t}{t}\right)^{\frac{1}{\theta}}  < +\infty .
		\end{align*}
		A measurable function $f:\Omega \rightarrow \mathbb{R}$ is said to belong to the Lorentz space $\mathcal{L}^{(p,\infty)} \left(\Omega\right)$ if 
		\begin{align*}
			\left\lVert f \right\rVert_{\mathcal{L}^{(p,\infty)}\left(\Omega\right)}:= \sup\limits_{t >0} t^{\frac{1}{p}}f^{\ast}_{\Omega}\left(t\right) < +\infty .
		\end{align*}
		Here $f_{\Omega}^{\ast}:[0, +\infty) \rightarrow [0, \infty)$ is the nonincreasing rearrangement function of $f$ over $\Omega,$ defined as 
		\begin{align*}
			f_{\Omega}^{\ast}\left( t \right):= \inf\left\lbrace s \geq 0: \left\lvert \left\lbrace x \in \Omega: \left\lvert f\left(x\right)\right\rvert > s\right\rbrace \right\rvert \leq t \right\rbrace. 
		\end{align*}
		
	\end{definition}
	The functions $\left\lVert \cdot  \right\rVert_{\mathcal{L}^{(p,\theta)}\left(\Omega\right)}$ and $\left\lVert \cdot  \right\rVert_{\mathcal{L}^{(p,\infty)}\left(\Omega;\varLambda^{k}\right)}$  in general defines only a quasinorm on the corresponding Lorentz spaces, which is not a norm. However, when $1 < p < \infty$ and $1 \leq \theta \leq \infty,$ the Lorentz quasinorm is equivalent to a norm which makes $\mathcal{L}^{(p,\theta)} \left(\Omega\right)$ into a Banach space ( see \cite{Bennet_Sharpley_InterpolationOperators} ). We would work only with these cases and hence would pretend that the quasinorm is actually a norm. For different properties of Lorentz spaces, see \cite{Bennet_Sharpley_InterpolationOperators}. The important point about the Lorentz spaces is that they are interpolation spaces, i.e. they can be obtained via real interpolation from the $L^{p}$ spaces. On the other hand, the Morrey spaces are not interpolation spaces. Roughly speaking, Morrey-Lorentz spaces are simply, Morrey type spaces based on Lorentz spaces, instead of the Lebesgue spaces for standard Morrey spaces. 
	\begin{definition}[Morrey-Lorentz spaces]
		Let  $1 <p < \infty, 1 \leq \theta < \infty, 0 \leq \mu < n $ be real numbers. For any measurable function $f:\Omega \rightarrow \mathbb{R},$ the Morrey-Lorentz quasinorm is defined as 
		\begin{align*}
			\left\lVert f\right\rVert_{\mathrm{L}^{(p,\theta)}_{\mu}\left(\Omega\right)}&:=\sup\limits_{\substack{z \in \Omega,\\ 0 < \rho \leq \operatorname{diam} \Omega}} \rho^{-\frac{\mu}{p}} \left(\int_{0}^{\infty}\left[ t^{\frac{1}{p}}f^{\ast}_{\Omega_{\left(\rho, z\right)}}\left(t\right) \right]^{\theta} \frac{\mathrm{d}t}{t}\right)^{\frac{1}{\theta}} \intertext{ and }
			\left\lVert f\right\rVert_{\mathrm{L}^{(p,\infty)}_{\mu}\left(\Omega\right)} &:= 
			\sup\limits_{\substack{z \in \Omega,\\ 0 < \rho \leq \operatorname{diam} \Omega}} \rho^{-\frac{\mu}{p}}\sup\limits_{t>0}\left[ t^{\frac{1}{p}}f^{\ast}_{\Omega_{\left(\rho, z\right)}}\left(t\right) \right] .
		\end{align*}
		We define the Morrey-Lorentz space 
		\begin{align*}
			\mathrm{L}^{(p, \theta)}_{\mu}\left(\Omega\right) := \left\lbrace  f:\Omega \rightarrow \mathbb{R} \text{ measurable }: \left\lVert f\right\rVert_{\mathrm{L}^{(p,\theta)}_{\mu}\left(\Omega\right)} < +\infty \right\rbrace.  
		\end{align*}
	\end{definition}
	Note that the norms can alternatively be expressed in terms of the Lorentz quasinorms as 
	\begin{align*}
		\left\lVert f\right\rVert_{\mathrm{L}^{(p,\theta)}_{\mu}\left(\Omega\right)} := \sup\limits_{\substack{z \in \Omega,\\ 0 < \rho \leq \operatorname{diam} \Omega}} \rho^{-\frac{\mu}{p}}\left\lVert f \right\rVert_{\mathcal{L}^{(p,\theta)}\left(\Omega_{\left(\rho,z\right)}\right)}. 
	\end{align*}
	The definition is extended componentwise in the obvious manner for $X$-valued functions, when $X$ is a finite dimensional real vector space. To avoid burdening our notations even more, we would often suppress the target space $X.$ 
	\subsubsection{H\"{o}lder inequality in Morrey-Lorentz spaces}	
	Now we record a H\"{o}lder inequality for Morrey-Lorentz spaces, which follows easily from the H\"{o}lder inequality for Lorentz spaces, due to Hunt \cite{Hunt_LorentzSpaces}. 
	\begin{theorem}
		Let $1 < p_{1}, p_{2} < \infty,$  $1 \leq \theta_{1}, \theta_{2} \leq \infty$ and $ 0 \leq \mu_{1}, \mu_{2} < n.$ Then for any $f \in \mathrm{L}_{\mu_{1}}^{\left(p_{1}, \theta_{1}\right)}$ and any $g \in \mathrm{L}_{\mu_{2}}^{\left(p_{2}, \theta_{2}\right)},$ we have $fg  \in \mathrm{L}_{\mu}^{\left(p, \theta\right)}$ and we have the estimate 
		\begin{align*}
			\left\lVert fg \right\rVert_{\mathrm{L}^{\left(p, \theta\right)}_{\mu}} \leq \left\lVert f \right\rVert_{\mathrm{L}^{\left(p_{1}, \theta_{1}\right)}_{\mu_{1}}}\left\lVert g \right\rVert_{\mathrm{L}^{\left(p_{2}, \theta_{2}\right)}_{\mu_{2}}},
		\end{align*}
		where 
		\begin{align*}
			\frac{1}{p} = 	\frac{1}{p_{1}} + 	\frac{1}{p_{2}},   \quad 	\frac{1}{\theta}	= \frac{1}{\theta_{1}} + 	\frac{1}{\theta_{2}} \quad \text{ and } \quad \frac{\mu}{p} = 	\frac{\mu_{1}}{p_{1}} + 	\frac{\mu_{2}}{p_{2}}. 
		\end{align*}
	\end{theorem}
	We would often use an important fact that the scaling of Lorentz norms is independent of the second exponent $\theta.$ 
	\begin{proposition}
		Let $r>0$ and let $1 < p < \infty$ and $1 \leq \theta \leq \infty.$ Let $U \subset \mathbb{R}^{n}$ be an open set and set $rU:= \left\lbrace rx: x \in U\right\rbrace.$ For any $u \in \mathcal{L}^{\left(p, \theta\right)}\left( rU\right).$  Then   
		\begin{align*}
			\left\lVert u\left( rx\right) \right\rVert_{\mathcal{L}^{\left(p, \theta\right)}\left(U\right)} = r^{-\frac{n}{p}} 	\left\lVert u \right\rVert_{\mathcal{L}^{\left(p, \theta\right)}\left(rU\right)}.  
		\end{align*}
		In particular, for any open, bounded subset $A \subset \mathbb{R}^{n},$ there exists a constant $C>0$ such that 
		\begin{align*}
			\left\lVert \mathbbm{1}_{A}\right\rVert_{\mathcal{L}^{\left(p, \theta\right)}\left( \mathbb{R}^{n}\right)} = C \left\lvert A \right\rvert^{\frac{1}{p}}. 
		\end{align*} 
	\end{proposition}
	\begin{remark}
		The constant $C$ depends on $\theta,$ but the power of $\left\lvert A \right\rvert$ does not. 
	\end{remark}
	\subsubsection{Sobolev spaces of Morrey-Lorentz type} We would also need Sobolev type spaces based on Morrey-Lorentz spaces.   
	\begin{definition}
		Let $l \geq 1$ be an integer and let  $1 <p < \infty, 1 \leq \theta < \infty, 0 \leq \mu < n $ be real numbers. The Morrey-Lorentz Sobolev spaces of order $l$ on $\Omega$ is defined as 
		\begin{align*}
			\mathsf{W}^{l}&\mathrm{L}^{(p,\theta)}_{\mu}\left( \Omega\right) := \left\lbrace u  \in \mathrm{L}^{(p,\theta)}_{\mu}\left( \Omega\right): D^{\alpha}u\in \mathrm{L}^{(p,\theta)}_{\mu}\left( \Omega\right) \text{ for all } 0 \leq \left\lvert \alpha \right\rvert \leq l  \right\rbrace,
		\end{align*}
		where $D^{\alpha}u$ denotes the $\alpha$-th weak derivative of $u.$ The space is equipped with the quasinorm 
		\begin{align*}
			\left\lVert u\right\rVert_{\mathsf{W}^{l}\mathrm{L}^{(p,\theta)}_{\mu}\left(\Omega\right)}:= \sum\limits_{0 \leq \left\lvert \alpha \right\rvert \leq l} 	\left\lVert D^{\alpha}u\right\rVert_{\mathsf{W}^{l}\mathrm{L}^{(p,\theta)}_{\mu}\left(\Omega\right)}. 
		\end{align*} 
	\end{definition}
	Once again, for any  finite dimensional real vector space $X$, $X$-valued Morrey-Lorentz Sobolev spaces  are defined componentswise and we would often write $\mathsf{W}^{l}\mathrm{L}^{(p, \theta)}_{\mu}\left(\Omega\right)$ in place  of $\mathsf{W}^{l}\mathrm{L}^{(p, \theta)}_{\mu}\left(\Omega;X\right).$ Also, we set $\mathsf{W}^{0}\mathrm{L}^{(p, \theta)}_{\mu}\left(\Omega\right) := \mathrm{L}^{(p, \theta)}_{\mu}\left(\Omega\right).$
	\begin{remark} Note that  
		\begin{enumerate}[(i)]  
			\item If $\mu=0$ and $1 < p=\theta < \infty,$ then these spaces are just Sobolev spaces based on Lebesgue spaces $L^{p},$ i.e. 
			\begin{align*}
				\mathrm{L}^{(p,p)}_{0}\left(\Omega\right) &\simeq L^{p}\left(\Omega\right) \quad \text{ with equivalent norms} &&\text{ for } l=0, \\
				\mathsf{W}^{l}\mathrm{L}^{(p,p)}_{0}\left(\Omega\right) &\simeq W^{l,p}\left(\Omega\right) \quad \text{ with equivalent norms} &&\text{ for } l\geq 1.
			\end{align*}
			\item If $\mu=0,$ $1 < p < \infty$ and $1 \leq \theta \leq \infty,$ then these spaces are the usual Lorentz-Sobolev spaces based on $\mathcal{L}^{(p,\theta)},$ i.e. 
			\begin{align*}
				\mathrm{L}^{(p,\theta)}_{0}\left(\Omega\right) &\simeq \mathcal{L}^{(p,\theta)}\left(\Omega\right) \quad \text{ with equivalent norms} &&\text{ for } l=0, \\
				\mathsf{W}^{l}\mathrm{L}^{(p,\theta)}_{0}\left(\Omega\right) &\simeq \mathsf{W}^{l}\mathcal{L}^{(p,\theta)}\left(\Omega\right) \quad \text{ with equivalent norms} &&\text{ for } l \geq 1.
			\end{align*}
			Also, when $\theta=\infty,$ they becomes Sobolev spaces based on Marcinkiewicz spaces or the weak Lebesgue spaces $L^{p}_{\text{w}}\left(\Omega\right),$ i.e.  
			\begin{align*}
				\mathrm{L}^{(p,\infty)}_{0}\left(\Omega\right) &\simeq L^{p}_{\text{w}}\left(\Omega\right) \quad \text{ with equivalent norms}&&\text{ for } l=0, \\
				\mathsf{W}^{l}\mathrm{L}^{(p,\infty)}_{0}\left(\Omega\right) &\simeq 	\mathsf{W}^{l}L^{p}_{\text{w}}\left(\Omega\right) \quad \text{ with equivalent norms}&&\text{ for } l \geq 1.
			\end{align*}
			\item If $ 0 < \mu <n$ and $1 < p=\theta < \infty,$ then these spaces are the usual Morrey-Sobolev spaces $\mathrm{L}^{p}_{\mu},$ i.e. 
			\begin{align*}
				\mathrm{L}^{(p,p)}_{\mu}\left(\Omega\right) &\simeq \mathrm{L}^{p}_{\mu}\left(\Omega\right) \quad \text{ with equivalent norms} &&\text{ for } l=0, \\
				\mathsf{W}^{l}\mathrm{L}^{(p,p)}_{\mu}\left(\Omega\right) &\simeq \mathsf{W}^{l}\mathrm{L}^{p}_{\mu}\left(\Omega\right) \quad \text{ with equivalent norms}&&\text{ for } l \geq 1.
			\end{align*}
			\item If $ 0 < \mu <n$, $1 < p < \infty$ and $\theta=\infty,$ then these spaces are the Sobolev spaces based on the so-called weak-Morrey spaces $\mathrm{L}_{\mu, \text{w}}^{p},$ i.e. 
			\begin{align*}
				\mathrm{L}^{(p,\infty)}_{\mu}\left(\Omega\right) &\simeq \mathrm{L}_{\mu, \text{w}}^{p}\left(\Omega\right) \quad \text{ with equivalent norms} &&\text{ for } l=0, \\
				\mathsf{W}^{l}\mathrm{L}^{(p,\infty)}_{\mu}\left(\Omega\right) &\simeq 	\mathsf{W}^{l}\mathrm{L}_{\mu, \text{w}}^{p}\left(\Omega\right) \text{ with equivalent norms} &&\text{ for } l \geq 1.
			\end{align*}
		\end{enumerate}
		
	\end{remark}
	\section{Sobolev and Poincar\'{e}-Sobolev inequalities}\label{Sobolev inequalities in Morrey-Lorentz}
	\subsection{Morrey-Lorentz Sobolev embeddings}
	We start with a result about extensions. Unfortunately, it would be too much of a digression to give a full proof here and it is difficult to find a reference in our particular setting. So we just sketch the basic ideas.  
	\begin{theorem}
		Let $1 < p <\infty,$ $0\leq \mu < n $  and $1 \leq \theta \leq \infty.$ Let $\Omega \subset \mathbb{R}^{n}$ be either a half ball or any open, bounded, smooth subset. For any integer $l \geq 0,$ there exists a linear bounded extension operator from $\mathsf{W}^{l}\mathrm{L}^{(p,\theta)}_{\mu}\left( \Omega\right) $ to $\mathsf{W}^{l}\mathrm{L}^{(p,\theta)}_{\mu}\left( \mathbb{R}^{n}\right).$ 
	\end{theorem}
	\begin{proof}
		We just sketch the ideas. Extension operator from $W^{l,p}\left( \Omega\right)$ to $W^{l,p}\left( \mathbb{R}^{n}\right)$ is standard and can be done in a number of different ways ( see \cite{Hestenes_ExtensionbyReflection},\cite{Calderon_SobolevSpaces}, \cite{Burenkov_SobolevSpacesDomains} ). By interpolation, these extend to bounded linear operators from $\mathsf{W}^{l}\mathcal{L}^{(p,\theta)}\left( \Omega\right) $ to $\mathsf{W}^{l}\mathcal{L}^{(p,\theta)}\left( \mathbb{R}^{n}\right).$ So one only needs to check that these operators preserve the Morrey-Sobolev type spaces as well. In our setting, where the domain is nice, this can be done. For pure Morrey-Sobolev cases, such results are proved in  \cite{Koskela_et_al_MorreySObolevExtension}, \cite{LambertiViolo_extensionMorreySobolev}, \cite{Lamberti_et_al_BurenkovExtensionMorreySobolev}, under far less regularity assumptions on the domain. 
	\end{proof}
	
	\begin{theorem}\label{Morrey-Lorentz subcritical}
		Let $1 < p <\infty,$ $0\leq \mu < n -p$  and $1 \leq \theta \leq \infty.$ Let $\Omega \subset \mathbb{R}^{n}$ be either a half ball or any open, bounded and $C^{1}$ subset. Then we have 
		\begin{align*}
			\mathsf{W}^{1}\mathrm{L}^{(p,\theta)}_{\mu}\left( \Omega\right) &\hookrightarrow 	\mathrm{L}^{\left( \frac{\left(n-\mu\right)p}{n - \mu -p},\frac{\left(n-\mu\right)\theta}{n-\mu -p}\right)}_{\mu}\left( \Omega\right) &&\text{ if } \mu >0, \\
			\mathsf{W}^{1}\mathcal{L}^{(p,\theta)}\left( \Omega\right) &\hookrightarrow 	\mathcal{L}^{\left( \frac{np}{n - p},\theta \right)}\left( \Omega\right) &&\text{ if } \mu = 0.
		\end{align*}
		
	\end{theorem}
	\begin{proof}
		Use the extension operator to extend any $u \in \mathsf{W}^{1}\mathrm{L}^{(p,\theta)}_{\mu}\left( \Omega\right)$	to a function $\tilde{u} \in \mathsf{W}^{1}\mathrm{L}^{(p,\theta)}_{\mu}\left( \mathbb{R}^{n}\right).$ Now the results follows from the boundedness of the fractional integral operators. In particular, we have 
		\begin{align*}
			I_{1}: \mathrm{L}^{(p,\theta)}_{\mu}\left(\mathbb{R}^{n}\right) \rightarrow \mathrm{L}^{\left( \frac{\left(n-\mu\right)p}{n - \mu -p},\frac{\left(n-\mu\right)\theta}{n-\mu -p}\right)}_{\mu}\left(\mathbb{R}^{n}\right) \qquad \text{ is bounded}.
		\end{align*} See Remark after Proposition 3 in \cite{Hatano_MorreyLorentz}. For the setting of Morrey spaces, i.e. the case $\theta =p,$ this result is due to Adams \cite{Adams_RieszPotentials}. Pure Lorentz case is due to Tartar \cite{Tartar_Lorentzembeddings}, Peetre \cite{Peetre_Lorentzembedding}.   
	\end{proof}
	As an immediate corollary, we record the following easy result. 
	\begin{theorem}
		Let  $1 < p <\infty,$ $\max\left\lbrace 0, n-p \right\rbrace \leq \mu < n$  and $1 \leq \theta \leq \infty.$ Let $\Omega \subset \mathbb{R}^{n}$ be either a half ball or any open, bounded and $C^{1}$ subset. Then we have 
		\begin{align*}
			\mathsf{W}^{1}\mathrm{L}^{(p,\theta)}_{\mu}\left( \Omega\right) \hookrightarrow 	L^{r}\left( \Omega\right)
		\end{align*}
		for every $1 \leq r < \infty.$	
	\end{theorem}
	
	\subsection{Poincar\'{e}-Sobolev inequalities} 
	\subsubsection{Poincar\'{e}-Sobolev inequalities for the gradient}
	We first record a simple compact embedding result. 
	\begin{proposition}[Compactness of embedding]
		Let $R>0$ and let $1 < p < \infty$ and $1 \leq \theta \leq \infty.$ Then the inclusion map from $\mathsf{W}^{1}\mathcal{L}^{\left(p, \theta\right)}\left(B_{R}^{+}\right)$ to $\mathcal{L}^{\left(p, \theta\right)}\left(B_{R}^{+}\right)$ is compact. 
	\end{proposition}
	\begin{proof}
		Choose $1 \leq  q < \min \left\lbrace p, n \right\rbrace$ and $\varepsilon >0$ such that $p < \frac{nq}{n-q} - \varepsilon.$ Now we have the continuous inclusions 
		\begin{align*}
			\mathsf{W}^{1}\mathcal{L}^{\left(p, \theta\right)}\left(B_{R}^{+}\right) \hookrightarrow W^{1,q}\left(B_{R}^{+}\right) \hookrightarrow 
			L^{\frac{nq}{n-q} - \varepsilon}\left(B_{R}^{+}\right) \hookrightarrow \mathcal{L}^{\left(p, \theta\right)}\left(B_{R}^{+}\right). 
		\end{align*}
		The claimed result follows as the middle inclusion is compact. 
	\end{proof}
	This compactness coupled with a simple contradiction argument proves the following two Poincar\'{e} inequalities. 
	\begin{proposition}[Poincar\'{e} inequality with zero mean in half balls for Lorentz spaces]\label{Poincare_zero_mean_Lorentz}
		Let $R>0$ and let $1 < p < \infty$ and $1 \leq \theta \leq \infty.$ Then for any $u \in \mathsf{W}^{1}\mathcal{L}^{\left(p, \theta\right)}\left(B_{R}^{+}\right)$ such that $ \fint_{B^{+}_{R}}u = 0,$ there exists a constant $C>0$ such that 
		\begin{align*}
			\left\lVert u \right\rVert_{\mathcal{L}^{\left(p, \theta\right)}\left(B_{R}^{+}\right)} \leq C R	\left\lVert \nabla u \right\rVert_{\mathcal{L}^{\left(p, \theta\right)}\left(B_{R}^{+}; \mathbb{R}^{n}\right)}.  
		\end{align*}
	\end{proposition}
	\begin{proposition}[Poincar\'{e} inequality in half balls for Lorentz spaces]\label{Poincare_zero_Dirichlet_Lorentz}
		Let $R>0$ and let $1 < p < \infty$ and $1 \leq \theta \leq \infty.$ Then for any $u \in \mathsf{W}^{1}\mathcal{L}^{\left(p, \theta\right)}\left(B_{R}^{+}\right)$ such that $u \equiv 0$ on $\Gamma_{R},$ there exists a constant $C>0$ such that 
		\begin{align*}
			\left\lVert u \right\rVert_{\mathcal{L}^{\left(p, \theta\right)}\left(B_{R}^{+}\right)} \leq C R	\left\lVert \nabla u \right\rVert_{\mathcal{L}^{\left(p, \theta\right)}\left(B_{R}^{+}; \mathbb{R}^{n}\right)}.  
		\end{align*}
	\end{proposition}
	These two Poincar\'{e} inequalities and Theorem \ref{Morrey-Lorentz subcritical} implies the following. 
	\begin{proposition}[Poincar\'{e}-Sobolev inequality in half balls for Lorentz spaces]\label{Poincare-Sobolev_Lorentz}
		Let $R>0$ and let $1 < p < n$ and $1 \leq \theta \leq \infty.$ Then for any $u \in \mathsf{W}^{1}\mathcal{L}^{\left(p, \theta\right)}\left(B_{R}^{+}\right)$ such that either $u \equiv 0$ on $\Gamma_{R},$ or $ \fint_{B^{+}_{R}}u = 0,$  there exists a constant $C>0$ such that 
		\begin{align*}
			\left\lVert u \right\rVert_{\mathcal{L}^{\left(\frac{np}{n-p}, \theta\right)}\left(B_{R}^{+}\right)} \leq C\left\lVert \nabla u \right\rVert_{\mathcal{L}^{\left(p, \theta\right)}\left(B_{R}^{+}; \mathbb{R}^{n}\right)}.  
		\end{align*}
	\end{proposition} 
	\subsubsection{Poincar\'{e}-Sobolev inequalities for the Hessian}
	\begin{lemma}\label{Hessian Poincare Sobolev}
		Let $R>0,$ $1 <p < n$ and $1 \leq \theta \leq \infty$. Then for any $u \in W^{1,2}\left( \mathcal{U}_{R}; \varLambda^{k}\right)$ satisfying 
		\begin{align*}
			\text{ either }\quad 	e_{n}\wedge u &=0  \qquad\text{ or } \qquad e_{n}\lrcorner u = 0 &&\text{ on } \Gamma_{3R/4},
		\end{align*} 
		and for any $0 < \rho \leq 3R/4,$ there exists $\bar{u}^{\rho} \in W^{1,2}\left( \mathcal{U}_{R}; \varLambda^{k}\right)$ such that  
		\begin{align*}
			D^{2}u &=D^{2}\bar{u}^{\rho} &&\text{ in } \mathcal{U}_{R}. 
		\end{align*} and there exists a constant $C = C \left( n, k, N, p, \theta \right) >0$ such that  
		\begin{align*}
			\frac{1}{\rho}\left\lVert \bar{u}^{\rho} \right\rVert_{\mathcal{L}^{(\frac{np}{n-p},\theta)}\left( B^{+}_{\rho}\right)} + \left\lVert \nabla \bar{u}^{\rho} \right\rVert_{\mathcal{L}^{(\frac{np}{n-p},\theta)}\left( B^{+}_{\rho}\right)} &\leq C\left\lVert D^{2}\bar{u}^{\rho} \right\rVert_{\mathcal{L}^{(p,\theta)}\left( B^{+}_{\rho}\right)}, 
		\end{align*}
		whenever $D^{2}u \in \mathcal{L}^{\left(p, \theta\right)}\left( B^{+}_{\rho} ; \varLambda^{k}\otimes \mathbb{R}^{n\times n}\right).$ 
		
	\end{lemma}
	\begin{proof}
		We first prove the case $e_{n}\wedge u = 0$ on $\Gamma_{3R/4}.$ By a simple scaling argument, we can assume $R=1$ and $ 0 < \rho \leq 3/4.$ Let us define 
		\begin{align*}
			\bar{u}^{\rho}_{I, j}\left( x\right) := \left\lbrace \begin{aligned}
				&u_{I, j}\left( x\right) - \left( \fint_{B^{+}_{\rho}} \frac{\partial 	u_{I, j}}{\partial x_{n}}\right)x_{n} &&\text{ if } n \notin I, \\
				&\begin{multlined}[t]
					u_{I, j}\left( x\right) - \left( \fint_{B^{+}_{\rho}}u_{I, j}\right) - \left\langle  x, \left( \fint_{B^{+}_{\rho}} \nabla u_{I, j}\right)\right\rangle \\+  \left\langle  \left( \fint_{B^{+}_{\rho}} x \right) , \left( \fint_{B^{+}_{\rho}} \nabla u_{I, j}\right)\right\rangle
				\end{multlined} &&\text{ if } n \in I,
			\end{aligned}\right. 
		\end{align*}
		for all $1 \leq j \leq N,$ where $\fint_{B^{+}_{\rho}} x$ denotes the constant vector in $\mathbb{R}^{n}$ formed by the components 
		$$ \left( \fint_{B^{+}_{\rho}} x\right)_{i} := \fint_{B^{+}_{\rho}} x_{i}\ \mathrm{d}x \qquad \text{ for } 1 \leq i \leq n.$$ From now on, every statement below is assumed to hold for every $1 \leq j \leq N.$ Now note that since $	e_{n}\wedge u =0$  on $ \Gamma_{\rho},$ we have $u_{I, j} \equiv 0$ on $ \Gamma_{\rho}$ if $ n \notin I$ and consequently, we also have 
		\begin{align*}
			\frac{\partial u_{I, j}}{\partial x_{l}} &\equiv 0 &&\text{ on } \Gamma_{\rho} \quad \text{ if } n \notin I, 1 \leq l \leq n-1.  
		\end{align*}
		Now it is easy to check that this implies, by our construction, that every component of $\bar{u}_{\rho}$ and all its first order derivatives either vanish on $\Gamma_{\rho}$ or has zero integral average on $B_{\rho}^{+}.$ The desired estimate easily follow from this by using the Propositions \ref{Poincare_zero_Dirichlet_Lorentz}, \ref{Poincare_zero_mean_Lorentz} and \ref{Poincare-Sobolev_Lorentz}, as appropriate. Since we also have $D^{2}u =D^{2}\bar{u}^{\rho} $ in $\mathcal{U},$ this completes the proof. For the case $e_{n}\lrcorner u = 0$ on $\Gamma_{3R/4},$ we interchange the cases $n \in I$ and $n \notin I$ in the definition of $\bar{u}^{\rho}_{I, j}$ and argue similarly. 
	\end{proof} 	
	\section{Crucial estimates}\label{crucial estimates}	
	\subsection{Lorentz estimates}
	\begin{theorem}\label{BoundaryLorentzzeroBVP}
		Let $l \geq 0$ be an integer and let  $1 <p < \infty, 1 \leq \theta \leq \infty.$ Let $\Omega \subset \mathbb{R}^{n}$ be a open, bounded, \textbf{contractible} subset such that $\partial\Omega$ is of class $C^{l+2}.$ Let $\bar{A} \in \operatorname{Hom}\left(\varLambda^{k+1}\right)$ and $B \in  \operatorname{Hom}\left(\varLambda^{k}\right)$ satisfy the Legendre condition. Then for any $f \in \mathsf{W}^{l}\mathcal{L}^{\left(p, \theta\right)}(\Omega, \varLambda^{k}),$ there exists a unique 
		solution $\omega \in  \mathsf{W}^{l+2}\mathcal{L}^{\left(p, \theta\right)}(\Omega, \varLambda^{k})$ to the following boundary value problem:
		\begin{equation}\label{bvp hodge elliptic system general full regularity Lorentz}
			\left\lbrace \begin{aligned}
				d^{\ast} ( \bar{A} d\omega )  + \left( \bar{B} \right)^{\intercal}d d^{\ast}\left( \bar{B} \omega \right)   &=  f  &&\text{ in } \Omega, \\
				\nu\wedge \omega &= 0 &&\text{  on } \partial\Omega. \\
				\nu\wedge d^{\ast} \left( \bar{B}\omega \right) &= 0 &&\text{ on } \partial\Omega, 
			\end{aligned} 
			\right. 
		\end{equation}
		which satisfies the estimate
		\begin{align*}
			\left\lVert \omega \right\rVert_{ \mathsf{W}^{l+2}\mathcal{L}^{\left(p, \theta\right)}(\Omega)} \leq c \left\lVert f\right\rVert_{ \mathsf{W}^{l}\mathcal{L}^{\left(p, \theta\right)}(\Omega)}.  
		\end{align*}
	\end{theorem}
	\begin{proof}
		First note that since $\Omega$ is contractible, by Theorem 16 and Remark 20 of \cite{Sil_linearregularity},  for any $f \in W^{l,r}\left( \Omega; \varLambda^{k}\right),$ there exists unique $\omega \in W^{l+2,r}\left(\Omega; \varLambda^{k}\right)$, solving \eqref{bvp hodge elliptic system general full regularity Lorentz} for any $1 < r< \infty$ and any integer $l \geq 0.$ Moreover, we have the estimate 
		\begin{align*}
			\left\lVert \omega \right\rVert_{W^{l+2,r}\left( \Omega; \varLambda^{k}\right)} \leq c \left( \left\lVert \omega \right\rVert_{W^{l,r}\left( \Omega; \varLambda^{k}\right)} + \left\lVert f \right\rVert_{W^{l,r}\left( \Omega; \varLambda^{k}\right)}\right). 
		\end{align*} Since the solution is unique, a simple contradiction and compactness argument implies that we in fact have the estimate 
		\begin{align*}
			\left\lVert \omega \right\rVert_{W^{l+2,r}\left( \Omega; \varLambda^{k}\right)} \leq c  \left\lVert f \right\rVert_{W^{l,r}\left( \Omega; \varLambda^{k}\right)}. 
		\end{align*}
		Note that this implies that the linear map $T$, defined by 
		\begin{align*}
			T\left(f\right): = \nabla^{2}\omega,
		\end{align*}
		where $\omega$ is the unique solution of \eqref{bvp hodge elliptic system general full regularity Lorentz} is a linear bounded operator from
		$L^{r}\left( \Omega; \varLambda^{k}\right)$ to $L^{r}\left( \Omega; \varLambda^{k}\otimes\mathbb{R}^{n\times n}\right)$ for any $1 < r < \infty.$ Now the general form of the Marcinkiewicz interpolation theorem ( see Theorem 4.13 in \cite{Bennet_Sharpley_InterpolationOperators} ) implies that this map extends as a bounded linear operator from $\mathcal{L}^{p,\theta}\left( \Omega; \varLambda^{k}\right)$ to $\mathcal{L}^{p, \theta}\left( \Omega; \varLambda^{k}\otimes\mathbb{R}^{n\times n}\right)$ for every $1 < p < \infty,$ $1 \leq \theta \leq \infty.$ The claimed result now follows easily. 
	\end{proof}
	\begin{theorem}\label{BoundaryLorentzzeroBVPnormal}
		Let $l \geq 0$ be an integer and let  $1 <p < \infty, 1 \leq \theta \leq \infty.$ Let $\Omega \subset \mathbb{R}^{n}$ be a open, bounded, \textbf{contractible} subset such that $\partial\Omega$ is of class $C^{l+2}.$ Let $\bar{A} \in \operatorname{Hom}\left(\varLambda^{k+1}\right)$ and $B \in  \operatorname{Hom}\left(\varLambda^{k}\right)$ satisfy the Legendre condition. Then for any $f \in \mathsf{W}^{l}\mathcal{L}^{\left(p, \theta\right)}(\Omega, \varLambda^{k}),$ there exists a unique 
		solution $\omega \in  \mathsf{W}^{l+2}\mathcal{L}^{\left(p, \theta\right)}(\Omega, \varLambda^{k})$ to the following boundary value problem:
		\begin{equation}\label{bvp hodge elliptic system general full regularity Lorentz normal}
			\left\lbrace \begin{aligned}
				\left( \bar{B}^{-1} \right)^{\intercal} d^{\ast} ( \bar{A} d\left( \bar{B}^{-1}\omega \right) )  + d d^{\ast} \omega    &=  f  &&\text{ in } \Omega, \\
				\nu\lrcorner \omega &= 0 &&\text{  on } \partial\Omega, \\
				\nu\lrcorner d^{\ast} \left( \bar{A} d\left( \bar{B}^{-1}\omega \right) \right) &= 0 &&\text{ on } \partial\Omega. 
			\end{aligned} 
			\right. 
		\end{equation}
		which satisfies the estimate
		\begin{align*}
			\left\lVert \omega \right\rVert_{ \mathsf{W}^{l+2}\mathcal{L}^{\left(p, \theta\right)}(\Omega)} \leq c \left\lVert f\right\rVert_{ \mathsf{W}^{l}\mathcal{L}^{\left(p, \theta\right)}(\Omega)}.  
		\end{align*}
	\end{theorem}
	\begin{proof}
		We set $u= \bar{B}^{-1}\omega$ to note that the desired estimate for $\omega$ is equivalent to deriving estimates for $u,$ where $u$ is the unique solution to the system 
		\begin{align*}
			\left\lbrace \begin{aligned}
				d^{\ast} ( \bar{A} du )  + \left( \bar{B} \right)^{\intercal} d d^{\ast} \left( \bar{B} u \right)      &=  \left( \bar{B} \right)^{\intercal}f  &&\text{ in } \Omega, \\
				\nu\lrcorner \left( \bar{B} u\right) &= 0 &&\text{  on } \partial\Omega, \\
				\nu\lrcorner d^{\ast} \left( \bar{A} du \right)  &= 0 &&\text{ on } \partial\Omega. 
			\end{aligned} 
			\right. 
		\end{align*}
		But the estimate for $u$ follows by interpolation from Theorem 17 in \cite{Sil_linearregularity}.  
	\end{proof}
	\subsection{Decay estimates}
	\begin{theorem}[Boundary Hessian decay estimates]\label{boundary hessian decay}
		Let $R>0,$ $1 <p  < q < \infty$ and $1 \leq \theta \leq \infty. $  Let $\bar{A} \in \operatorname{Hom}\left(\varLambda^{k+1}\right),$ $\bar{B} \in \operatorname{Hom}\left(\varLambda^{k}\right)$ satisfy the Legendre condition with constant $\gamma >0.$ Assume one of the following holds.\smallskip 
		
		\begin{enumerate}[(i)] 
			\item Let $\alpha \in W^{1,2}(\mathcal{U}_{R} ; \varLambda^{k})\cap \mathsf{W}^{2}\mathcal{L}^{(p,\theta)}\left( \mathcal{U}_{R}; \varLambda^{k}\right) $  satisfy $e_{n}\wedge \alpha = 0$ on $\Gamma_{3R/4}$ and for all $\psi \in W_{T}^{1,2}\left( \mathcal{U}_{R} ; \varLambda^{k}\right),$ we have 
			\begin{align}\label{bndry gradient decay estimate constant}
				\int_{\mathcal{U}_{R}} \langle \bar{A}d\alpha ; d\psi \rangle + \int_{\mathcal{U}_{R}} \langle d^{\ast}\left( \bar{B} \alpha \right) ; d^{\ast} \left( \bar{B}\psi\right) \rangle   = 0.
			\end{align}
			\item Let $\alpha \in W^{1,2}(\mathcal{U}_{R} ; \varLambda^{k})\cap \mathsf{W}^{2}\mathcal{L}^{(p,\theta)}\left( \mathcal{U}_{R}; \varLambda^{k}\right) $  satisfy $e_{n}\lrcorner \alpha = 0$ on $\Gamma_{3R/4}$ and for all $\psi \in W_{N}^{1,2}\left( \mathcal{U}_{R} ; \varLambda^{k}\right),$ we have 
			\begin{align}\label{bndry gradient decay estimate constant normal}
				\int_{\mathcal{U}_{R}} \langle \bar{A}d\left( \bar{B}^{-1} \alpha \right) ; d\left( \bar{B}^{-1} \psi \right) \rangle + \int_{\mathcal{U}_{R}} \langle d^{\ast}\alpha ; d^{\ast} \psi \rangle   = 0.  
			\end{align}
		\end{enumerate}\noindent 
		Then we have 
		$D^{2}\alpha \in \mathcal{L}^{(q, \theta)}\left( B^{+}_{R/2};\varLambda^{k}\otimes\mathbb{R}^{n \times n}\right) $ and there exists a constant $C = C\left( p, q, \theta,  \gamma, k, n, N \right) >0$ such that we have the estimate
		\begin{equation}\label{boundary Hessian higher integrability}
			R^{\frac{n}{p}- \frac{n}{q}}\left\lVert D^{2}\alpha \right\rVert_{\mathcal{L}^{(q, \theta)}\left( B^{+}_{R/2}\right)} \leq C \left\lVert D^{2}\alpha \right\rVert_{\mathcal{L}^{(p, \theta)}\left( B^{+}_{3R/4}\right)}.
		\end{equation}
	\end{theorem}
	\begin{proof}
		We first show $(i)$. By scale invariance of \eqref{boundary Hessian higher integrability}, we can assume $R=1.$
		First assume $q <n.$ Since $p <q,$ there exists $m \in \mathbb{N}$ such that 
		$$ \frac{1}{p} - \frac{m}{n} \leq \frac{1}{q} < \frac{1}{p} - \frac{m-1}{n}. $$
		Now, for every $1 \leq j \leq m+1,$ define the radii $r_{j}$ and the exponents $q_{j}$ by 
		$$ r_{j}:= \frac{3}{4} - \frac{j-1}{4m}\qquad \text{ and } \qquad q_{j}:= \frac{np}{n - \left( j-1\right)p}.$$ 
		Now we claim that for every $1 \leq j \leq m,$ $\left\lvert D^2\alpha \right\rvert \in \mathcal{L}^{\left( q_{j}, \theta \right)}\left( B^{+}_{r_{j}}\right)$ and there exist constants $C_{j} >0,$ independent of $\alpha,$ such that we have the estimate 
		\begin{align*}
			\left\lVert D^{2}\alpha \right\rVert_{\mathcal{L}^{\left( q_{j+1}, \theta \right)}\left( B^{+}_{r_{j+1}}\right)}\leq C_{j} 	\left\lVert D^{2}\alpha \right\rVert_{\mathcal{L}^{\left( q_{j}, \theta \right)}\left( B^{+}_{r_{j}}\right)}. 
		\end{align*}
		The claim implies the result, as combining the estimates, we get 
		\begin{align*}
			\left\lVert D^{2}\alpha \right\rVert_{\mathcal{L}^{\left( q_{m+1}, \theta \right)}\left( B^{+}_{r_{m+1}}\right)} \leq \left( \prod_{j=1}^{m} C_{j} \right) \	\left\lVert D^{2}\alpha \right\rVert_{\mathcal{L}^{\left( q_{1}, \theta \right)}\left( B^{+}_{r_{1}}\right)}
		\end{align*}
		But this is our desired estimate as $q_{1}=p,$ $r_{1} = 3/4,$ $r_{m+1} = 1/2 $ and $q \leq q_{m+1}.$ We prove the claim by induction. Fix $1 \leq j \leq m$ and assume the claim holds for all $1 \leq l \leq j-1.$ Thus, we have $D^{2}\alpha \in \mathcal{L}^{\left( q_{j}, \theta \right)}\left( B^{+}_{r_{j}}\right).$ Thus, using $\rho = r_{j}$ in Lemma \eqref{Hessian Poincare Sobolev}, there exists $\bar{\alpha}^{j} \in W^{1,2}\left( \mathcal{U}; \varLambda^{k}\right)$ such that $D^{2}\alpha =D^{2}\bar{\alpha}^{j}$ in $\mathcal{U}$ and 
		\begin{align*}
			\frac{1}{r_{j}}\left\lVert \bar{\alpha}^{j} \right\rVert_{\mathcal{L}^{(q_{j+1},\theta)}\left( B^{+}_{r_{j}}\right)} + \left\lVert \nabla \bar{\alpha}^{j} \right\rVert_{\mathcal{L}^{(q_{j+1},\theta)}\left( B^{+}_{r_{j}}\right)} &\leq C\left\lVert D^{2}\bar{\alpha}^{j} \right\rVert_{\mathcal{L}^{(q_{j},\theta)}\left( B^{+}_{r_{j}}\right)} 
		\end{align*}
		for some constant $C >0.$ Now choose a scalar cut-off functions $\zeta_{j}: B_{1} \rightarrow \mathbb{R}$ such that 
		\begin{gather*}
			\zeta_{j} \in C_{c}^{\infty}\left( B_{r_{j}}\right), \quad 0 \leq \zeta_{j} \leq 1 \text{ in } B_{r_{j}}, \quad \zeta_{j} \equiv 1 \text{ in } B_{r_{j+1}}, \\  
			\left\lvert \nabla \zeta_{j}\right\rvert \leq \frac{C}{\left( r_{j} -r_{j+1}\right)}, \quad \left\lvert D^{2} \zeta_{j}\right\rvert \leq \frac{C}{\left( r_{j} -r_{j+1}\right)^{2}}. \end{gather*}  for some fixed constant $ C>0.$ 	Now a direct calculation shows that $\beta_{j}:= \zeta_{j}\bar{\alpha}^{j}$ is a weak solution of \eqref{bvp hodge elliptic system general full regularity Lorentz} with $\Omega := \mathcal{U}$ and 
		\begin{multline*}
			f:=d^{\ast} \left[ \bar{A}\left( d\zeta_{j}\wedge \bar{\alpha}^{j}\right)\right] + \left( \bar{B}\right)^{\intercal}d \left(d\zeta_{j}\lrcorner \bar{B}\bar{\alpha}^{j}\right) - d\zeta_{j}\lrcorner \left[ \bar{A}d\bar{\alpha}^{j} \right] \\ - \left( \bar{B}\right)^{\intercal} \left[ d\zeta_{j}\wedge d^{\ast}\left( \bar{B}\bar{\alpha}^{j}\right)\right] . 
		\end{multline*}
		Easy calculations imply that $ F \in \mathcal{L}^{\left( q_{j+1}, \theta \right)}\left( \mathcal{U}; \varLambda^{k}\right)$ along with the estimate 
		\begin{align*}
			\left\lVert F \right\rVert_{\mathcal{L}^{\left( q_{j+1}, \theta \right)}\left(  \mathcal{U}\right)}	\leq C_{j}\left\lVert D^{2}\bar{\alpha}^{j}  \right\rVert_{\mathcal{L}^{\left( q_{j}, \theta \right)}\left( B^{+}_{r_{j}}\right)} , 
		\end{align*}
		where the constant $C_{j}>0$ depends  on $j.$ Applying Theorem \ref{BoundaryLorentzzeroBVP}, we obtain the estimate 
		\begin{align*}
			\left\lVert \beta_{j}\right\rVert_{\mathsf{W}^{2}\mathcal{L}^{\left( q_{j+1}, \theta \right)} \left(\mathcal{U}\right)} \leq C \left\lVert F \right\rVert_{\mathcal{L}^{\left( q_{j+1}, \theta \right)} \left(\mathcal{U}\right)} \leq C_{j} \left\lVert D^{2}\bar{\alpha}^{j} \right\rVert_{\mathcal{L}^{\left( q_{j}, \theta \right)} \left(B^{+}_{r_{j}}\right)}  .
		\end{align*} 
		Now since $\zeta_{j} \equiv 1$ in $B^{+}_{r_{j+1}},$ we deduce 
		\begin{align*}
			\left\lVert D^{2}\bar{\alpha}^{j} \right\rVert_{\mathcal{L}^{\left( q_{j+1}, \theta \right)} \left(B^{+}_{r_{j+1}}\right)} &\leq  \left\lVert \beta_{j}\right\rVert_{\mathsf{W}^{2}\mathcal{L}^{\left( q_{j+1}, \theta \right)} \left(\mathcal{U}\right)} \leq C_{j} \left\lVert D^{2}\bar{\alpha}^{j} \right\rVert_{\mathcal{L}^{\left( q_{j}, \theta \right)} \left(B^{+}_{r_{j}}\right)} . 
		\end{align*}
		This proves the estimate in case $q<n.$ Other cases are easier. For $(ii),$ we just use Theorem \ref{BoundaryLorentzzeroBVPnormal} instead of Theorem \ref{BoundaryLorentzzeroBVP}. \end{proof}
	\subsection{Flattening the boundary}
	By obvious modifications in the proof of Lemma $4$ in  \cite{Sil_linearregularity}, we have the following. 
	\begin{lemma}\label{flattening and freezing coefficients}
		Let $l \geq 0$ be an integer and let  $1 <p < \infty, 1 \leq \theta \leq \infty, 0 \leq \mu < n .$  Let $\Omega \subset \mathbb{R}^{n}$ be a open, bounded subset such that $\partial\Omega$ is of class $C^{l+2}.$ Let $A \in C^{l+1}\left(\overline{\Omega}; \operatorname{Hom}\left(\varLambda^{k+1}\right)\right)$ and $B \in C^{l+2}\left(\overline{\Omega}; \operatorname{Hom}\left(\varLambda^{k}\right)\right)$ satisfy the Legendre condition. Suppose  $f \in \mathsf{W}^{l}\mathrm{L}^{(p,\theta)}_{\mu}\left( \Omega; \varLambda^{k}\right) \cap W^{l,2}\left( \Omega; \varLambda^{k}\right)$ and $\lambda \in \mathbb{R}.$ \smallskip 
		
		\noindent Let $\omega \in W_{T}^{1,2}(\Omega, \varLambda^{k}) \cap \mathsf{W}^{l+2}\mathcal{L}^{(p,\theta)}\left(\Omega; \varLambda^{k}\right) \cap W^{l+2,2}\left( \Omega; \varLambda^{k}\right)$ satisfy, 
		\begin{align}\label{bilinear form w22 regularity}
			\int_{\Omega}  \langle A (x)d\omega , d\phi \rangle   + \int_{\Omega}\langle d^{\ast} \left( B(x)\omega \right) , d^{\ast} \left( B(x)\phi \right) \rangle 
			+ \lambda  \int_{\Omega}\langle  B(x)\omega , \phi \rangle 		+\int_{\Omega} \langle  f, \phi \rangle = 0,    
		\end{align}
		for all $ \phi \in W_{T}^{1,2}(\Omega; \varLambda^{k}).$ Then for every $x_0\in\partial\Omega,$  there exist  
		\begin{enumerate}[(i)]
			\item a positive number $0 < R_{0}<1$ and a neighborhood $U$ of $x_0$ in $\mathbb{R}^n$ and  such that there exists $\Phi\in \operatorname{Diff}^{l+2}(\overline{B_{R_{0}}};\overline{U})$ with 
			\begin{align*}
				\Phi(0)=x_0, \quad D\Phi(0) \in \mathbb{SO}\left(n\right), \quad\Phi(B^{+}_{R_{0}})=\Omega\cap U, \quad  \Phi(\Gamma_{R_{0}})=\partial\Omega\cap U, 
			\end{align*}
			\item a scalar function $\zeta \in C_{c}^{\infty}\left( U\right)$ and constant matrices $\bar{A} \in \operatorname{Hom}\left( \varLambda^{k+1} \right)$ and $\bar{B}\in \operatorname{Hom}\left( \varLambda^{k} \right),$ both satisfying the Legendre condition, 
			\item vector-valued functions  $ \mathrm{P} ,\mathrm{Q} ,\mathrm{R}$  and  $\mathrm{S} $ with 
			\begin{align*}
				& \mathrm{P} \in C^{l}\left(\overline{B^{+}_{R_{0}}}; \operatorname{Hom}\left( \varLambda^{k}\right)\right), \mathrm{Q} \in C^{l}\left(\overline{B^{+}_{R_{0}}}; \operatorname{Hom}\left( \varLambda^{k}; \varLambda^{k}\otimes \mathbb{R}^{n}\right)\right), \\& \mathrm{R} \in C^{l}\left(\overline{B^{+}_{R_{0}}}; \operatorname{Hom}\left( \varLambda^{k}\otimes \mathbb{R}^{n}; \varLambda^{k}\right)\right) \text{ and } \mathrm{S} \in C^{l+1}\left(\overline{B^{+}_{R_{0}}}; \operatorname{Hom}\left( \varLambda^{k}\otimes \mathbb{R}^{n}\right)\right),
			\end{align*}
			depending only on $A$, $B$, $\Phi$, $\zeta$, $U$ and $R_{0},$ such that 
			\begin{align}\label{smallness of S}
				\left\lVert \mathrm{S} \right\rVert_{L^{\infty}\left(B_{r}^{+}\right)} \leq Cr \qquad \text{ for all } 0 < r \leq R_{0},  
			\end{align}
			\item $\widetilde{f} \in \mathsf{W}^{l}\mathrm{L}^{(p,\theta)}_{\mu}\left( B^{+}_{R_{0}}; \varLambda^{k}\right) \cap  W^{l,2}\left( B^{+}_{R_{0}}; \varLambda^{k}\right),$ with estimates on the $\mathsf{W}^{l}\mathrm{L}^{(p,\theta)}$ and $W^{l,2}$ norms by the corresponding norms of $f$, with the constants in the estimates depending only on $\Phi$, $\zeta$, $U$ and $R_{0}$,
		\end{enumerate}
		such that for all $\psi \in W_{T, flat}^{1,2}(B_{R_{0}}^{+} ; \Lambda^{k}),$ we have 
		\begin{align*}
			\int_{B_{R_{0}}^{+}} \langle \bar{A}(du) ; d\psi \rangle + \int_{B_{R_{0}}^{+}} &\langle d^{\ast} ( \bar{B} u ) ; d^{\ast} ( \bar{B} \psi )  \rangle 
			+  \int_{B_{R_{0}}^{+}} \langle \widetilde{f} + \mathrm{P}u + \mathrm{R}\nabla u ; \psi \rangle 
			\notag \\ & + \int_{B_{R_{0}}^{+}} \langle \mathrm{Q} u ; \nabla\psi \rangle   + \int_{B_{R_{0}}^{+}} \langle \mathrm{S} \nabla u , \nabla\psi \rangle = 0,
		\end{align*}
		where $u = \Phi^{\ast}(\zeta\omega) \in W_{T, flat}^{1,2}\left(B_{R_{0}}^{+} ; \varLambda^{k}\right) \cap \mathsf{W}^{l+2}\mathcal{L}^{(p,\theta)}\left(B_{R_{0}}^{+} ; \varLambda^{k}\right) \cap W^{l+2,2}\left(B_{R_{0}}^{+} ; \varLambda^{k}\right).$ 
	\end{lemma}
	\begin{remark}\label{flatteninglemmanormal} If $\omega \in W_{N}^{1,2}(\Omega, \varLambda^{k})$ satisfy, 
		\begin{multline}
			\int_{\Omega}  \langle A (x)d \left( B^{-1}(x) \omega \right) , d \left( B^{-1}(x) \phi \right) \rangle   + \int_{\Omega}\langle d^{\ast} \omega  , d^{\ast} \phi \rangle 
			+ \lambda  \int_{\Omega}\langle  \omega , B^{-1}(x)\phi \rangle \\
			+\int_{\Omega} \langle  f, B^{-1}(x)\phi \rangle   = 0,   
		\end{multline}
		for all $ \phi \in W_{N}^{1,2}(\Omega; \varLambda^{k})$ then analogous results hold, giving the existence $W,\theta, \Phi$ and constant matrices  $\bar{A}$ and $\bar{B}$, both satisfying 
		the Legendre condition such that $u = \Phi^{\ast}(\theta\omega) \in W_{N, flat}^{1,2}(B_{R}^{+} ; \varLambda^{k}) $  satisfies, 
		for all $\psi \in W_{N, flat}^{1,2}(B_{R}^{+} ; \varLambda^{k}),$
		\begin{align*}
			\int_{B_{R}^{+}} \langle \bar{A}(d \left( \bar{B}^{-1}u \right)) ; d\left( \bar{B}^{-1}\psi \right)\rangle + &\int_{B_{R}^{+}} \langle d^{\ast}  u ; d^{\ast} \psi   \rangle 
			+  \int_{B_{R}^{+}} \langle \widetilde{f} + \mathrm{P}u + \mathrm{R}\nabla u ; \psi \rangle 
			\notag \\ & + \int_{B_{R}^{+}} \langle \mathrm{Q} u ; \nabla\psi \rangle   + \int_{B_{R}^{+}} \langle \mathrm{S} \nabla u , \nabla\psi \rangle = 0,
		\end{align*}
		with the same conclusions for $\mathrm{P},\mathrm{Q},\mathrm{R},\mathrm{S}$ and $\widetilde{f}.$ 
	\end{remark}
	\subsection{Boundary estimates}
	\begin{lemma}\label{boundary estimate lemma}
		Let $l \geq 0$ be an integer and let  $1 <p < \infty, 1 \leq \theta \leq \infty, 0 \leq \mu < n .$  Let $\Omega \subset \mathbb{R}^{n}$ be a open, bounded subset such that $\partial\Omega$ is of class $C^{l+2}.$ Let $A \in C^{l+1}\left(\overline{\Omega}; \operatorname{Hom}\left(\varLambda^{k+1}\right)\right)$  and $B \in C^{l+2}\left(\overline{\Omega}; \operatorname{Hom}\left(\varLambda^{k+1}\right)\right)$ satisfy the Legendre condition. Suppose  $f \in \mathsf{W}^{l}\mathrm{L}^{(p,\theta)}_{\mu}\left( \Omega; \varLambda^{k}\right) \cap W^{l,2}\left( \Omega; \varLambda^{k}\right)$  and $\lambda \in \mathbb{R}.$ Let $\omega \in  W^{l+2,2}\left( \Omega; \varLambda^{k}\right) \cap \mathsf{W}^{l+2}\mathrm{L}_{\kappa}^{(p,\theta)}\left(\Omega; \varLambda^{k}\right)$ and for any $ 0 \leq \kappa < \mu <n,$ let us set  
		\begin{align*}
			\tilde{\kappa} := \min \left\lbrace \kappa + p, \mu \right\rbrace.  
		\end{align*}
		Assume one of the following holds. 
		\begin{enumerate}[(i)]
			\item Let $\omega \in W_{T}^{1,2}(\Omega, \varLambda^{k})$ and for all  $ \phi \in W_{T}^{1,2}(\Omega; \varLambda^{k}),$ we have, 
			\begin{multline}\label{bilinear form boundary estimate}
				\int_{\Omega}  \langle A (x)d\omega , d\phi \rangle   + \int_{\Omega}\langle d^{\ast} \left( B(x)\omega \right) , d^{\ast} \left( B(x)\phi \right) \rangle 
				+ \lambda  \int_{\Omega}\langle  B(x)\omega , \phi \rangle 	\\ +\int_{\Omega} \langle  f, \phi \rangle   = 0.    
			\end{multline}
			\item Let $\omega \in W_{B,N}^{1,2}(\Omega, \varLambda^{k})$ and for all $ \phi \in W_{B,N}^{1,2}(\Omega; \varLambda^{k}),$ we have, 
			\begin{multline}\label{bilinear form boundary estimate normal}
				\int_{\Omega}  \langle A (x)d\omega , d\phi \rangle   + \int_{\Omega}\langle d^{\ast} \left( B(x)\omega \right) , d^{\ast} \left( B(x)\phi \right) \rangle 
				+ \lambda  \int_{\Omega}\langle  B(x)\omega , \phi \rangle 	\\	+\int_{\Omega} \langle  f, \phi \rangle  = 0.    
			\end{multline}
		\end{enumerate}
		\noindent  Then for every $x_0\in\partial\Omega,$  there exist  
		there exists  $0 < R  < 1,$
		a neighborhood $U$ of $x_0$ in $\mathbb{R}^n$ and $\Phi\in \operatorname{Diff}^{l+2}(\overline{B_{R_{0}}};\overline{U})$, such that 
		\begin{align*}
			\Phi(0)=x_0,  D\Phi(0) \in \mathbb{SO}\left(n\right), \Phi(B^{+}_{R})=\Omega\cap U  \text{ and }  \Phi(\Gamma_R)=\partial\Omega\cap U, 
		\end{align*} and a constant $C = C \left(x_{0}, n, k, N, \gamma, \Omega, \lambda, \kappa, p, \theta, \mu  \right)>0,$ such that we have 
		\begin{align}\label{u boundary estimate}
			\left\lVert \omega \right\rVert_{\mathsf{W}^{l+2}\mathrm{L}^{(p,\theta)}_{\tilde{\kappa}}\left( \Phi \left( B^{+}_{R/2}\right)\right)}  \leq C \left( \left\lVert \omega \right\rVert_{\mathsf{W}^{l+2}\mathrm{L}_{\kappa}^{(p,\theta)}\left( \Omega\right)} + \left\lVert f \right\rVert_{\mathsf{W}^{l+2}\mathrm{L}^{(p,\theta)}_{\mu}\left( \Omega\right)} \right).
		\end{align}
	\end{lemma}
	\begin{proof}
		We begin with $(i).$ We prove only for $l =0$ as the result can be iterated. Using Lemma \ref{flattening and freezing coefficients}, for every $x_{0} \in \partial\Omega,$ there exists a positive number $0 < R_{0}<1,$ a neighborhood $U$ of $x_0$ in $\mathbb{R}^n$ and $\Phi\in \operatorname{Diff}^{2}(\overline{B_{R_{0}}};\overline{U})$ such that $u = \Phi^{\ast}(\zeta\omega) \in W_{T, flat}^{1,2}\left(B_{R_{0}}^{+} ; \varLambda^{k}\right) \cap \mathsf{W}^{2}\mathcal{L}^{(p,\theta)}\left(B_{R_{0}}^{+} ; \varLambda^{k}\right) \cap W^{2,2}\left(B_{R_{0}}^{+} ; \varLambda^{k}\right)$ satisfies, 
		\begin{align}\label{flattened equation}
			\int_{B_{R_{0}}^{+}} \langle \bar{A}(du) ; d\psi \rangle + \int_{B_{R_{0}}^{+}} &\langle d^{\ast} ( \bar{B} u ) ; d^{\ast} ( \bar{B} \psi )  \rangle 
			+  \int_{B_{R_{0}}^{+}} \langle \widetilde{f} + \mathrm{P}u + \mathrm{R}\nabla u ; \psi \rangle 
			\notag \\ & + \int_{B_{R_{0}}^{+}} \langle \mathrm{Q} u ; \nabla\psi \rangle   + \int_{B_{R_{0}}^{+}} \langle \mathrm{S} \nabla u , \nabla\psi \rangle = 0
		\end{align}
		for all $\psi \in W_{T, flat}^{1,2}(B_{R_{0}}^{+} ; \varLambda^{k}),$ where $\bar{A}, \bar{B}, \mathrm{P}, \mathrm{Q}, \mathrm{R}, \mathrm{S}, \tilde{f}, \Phi$ are as in Lemma \ref{flattening and freezing coefficients}. 
		
		Now  let $0 < R < R_{0}.$ We are going to choose $R$ later. The constants in all the estimates that we would derive from here onward may depend on $R_{0},$ but does not depend on $R.$ Since $\Phi$ is a diffeomorphism, it is enough to estimate $ 	\left\lVert u \right\rVert_{\mathsf{W}^{2}\mathrm{L}^{(p,\theta)}_{\tilde{\kappa}}\left( B^{+}_{R/2}\right)}.$ In view of \eqref{estimate for lower order terms}, it suffices to estimate   
		\begin{align*}
			\rho^{-\frac{\tilde{\kappa}}{p}} \left\lVert D^{2}u\right\rVert_{\mathcal{L}^{\left(p, \theta \right)}\left( B_{\rho}\left( y \right)\cap B^{+}_{R/2} \right)} \qquad \text{ for } y \in B^{+}_{R/2}.
		\end{align*}
		The estimate is trivial when $\rho$ has a lower bound, so we fix $0 < \sigma < R/32$ and show the estimate for $0 < \rho < \sigma/2.$ Let $y = \left( y', y_{n}\right) \in B^{+}_{R/2}.$ Denoting the point $\left( y', 0\right) \in \partial\mathbb{R}^{n}_{+}$ still by $y',$ we note that either  $y_{n} > \sigma,$ in which case $B_{\sigma}\left( y\right) \subset \subset B^{+}_{R},$ or we have  
		$0 \leq y_{n} \leq \sigma$ and then $B_{\sigma}\left( y\right) \cap B^{+}_{R} \subset B^{+}_{2\sigma}\left( y' \right) \subset B^{+}_{9R/16}\subset B^{+}_{3R/4}. $

		We only show the estimate for this last case, as the other is an interior estimate. By existence theory, there exists a $\beta\in W^{1,2}_{T}\left( \mathcal{U}_{2\sigma}\left(y'\right); \varLambda^{k}\right)$ such that 
		\begin{align}\label{frozen system for beta}
			\left\lbrace \begin{aligned}
				d^{\ast}\left( \bar{A}d\beta\right) + \bar{B}^{T}dd^{\ast}\left(\bar{B}\beta\right) &= g -\operatorname{div} G &&\text{ in } \mathcal{U}_{2\sigma}\left(y'\right), \\
				\nu \wedge \beta &=0  &&\text{ on } \partial\mathcal{U}_{2\sigma}\left(y'\right), \\
				\nu \wedge d^{\ast}\left(\bar{B}\beta\right) &=0  &&\text{ on } \partial\mathcal{U}_{2\sigma}\left(y'\right),
			\end{aligned}\right. 
		\end{align}
		where 
		\begin{align*}
			g := \widetilde{f} + \mathrm{P}u + \mathrm{R}\nabla u \qquad \text{ and } \qquad G := \mathrm{Q} u + S\nabla u . 
		\end{align*}
		Now it follows from Theorem \ref{BoundaryLorentzzeroBVP} that we have 
		\begin{align}\label{boundary Lorentz estimate for beta}
			\left\lVert D^{2}\beta\right\rVert_{\mathcal{L}^{^{(p,\theta)}}\left(\mathcal{U}_{2\sigma}\left(y'\right)\right)} \leq C 	\left\lVert g -\operatorname{div} G\right\rVert_{\mathcal{L}^{^{(p,\theta)}}\left(\mathcal{U}_{2\sigma}\left(y'\right)\right)} ,
		\end{align}
		where the constant $C$ is independent of $\sigma >0,$ as can be easily seen by scaling. Now we estimate the right hand side. We have 
		\begin{align*}
			\left\lVert \operatorname{div} \mathrm{Q}u \right\rVert_{\mathcal{L}^{^{(p,\theta)}}\left( \mathcal{U}_{2\sigma}\left(y'\right)\right)}  & \leq 	C\left( \left\lVert  u \right\rVert_{\mathcal{L}^{^{(p,\theta)}}\left( \mathcal{U}_{2\sigma}\left(y'\right)\right)} + \left\lVert \nabla u \right\rVert_{\mathcal{L}^{^{(p,\theta)}}\left( \mathcal{U}_{2\sigma}\left(y'\right)\right)} \right), \\
			\left\lVert g \right\rVert_{\mathcal{L}^{^{(p,\theta)}}\left( \mathcal{U}_{2\sigma}\left(y'\right)\right)} 
			&\leq  \begin{multlined}[t]
				C\left( \left\lVert  u\right\rVert_{\mathcal{L}^{^{(p,\theta)}}\left( \mathcal{U}_{2\sigma}\left(y'\right)\right)}   + \left\lVert  \nabla u\right\rVert_{\mathcal{L}^{^{(p,\theta)}}\left( \mathcal{U}_{2\sigma}\left(y'\right)\right)} \right)  \\ +	C\left\lVert \widetilde{f} \right\rVert_{\mathcal{L}^{^{(p,\theta)}}\left( \mathcal{U}_{2\sigma}\left(y'\right)\right)} 
			\end{multlined}
		\end{align*}
		and the last term is estimated as 
		\begin{multline*}
			\left\lVert \operatorname{div} \mathrm{S}\nabla u \right\rVert_{\mathcal{L}^{^{(p,\theta)}}\left( \mathcal{U}_{2\sigma}\left(y'\right)\right)}  \leq 	 \left\lVert S \right\rVert_{L^{\infty}\left(B^{+}_{R}\right)}\left\lVert  D^{2}u \right\rVert_{\mathcal{L}^{^{(p,\theta)}}\left( \mathcal{U}_{2\sigma}\left(y'\right)\right)} \\ + C\left\lVert \nabla u \right\rVert_{\mathcal{L}^{^{(p,\theta)}}\left( \mathcal{U}_{2\sigma}\left(y'\right)\right)}. 
		\end{multline*}
		Combining these estimates with \eqref{boundary Lorentz estimate for beta}, we deduce 
		\begin{align}
			\left\lVert D^{2}\beta\right\rVert_{\mathcal{L}^{^{(p,\theta)}}\left(\mathcal{U}_{2\sigma}\left(y'\right)\right)} \leq \begin{multlined}[t]
				C\left\lVert S \right\rVert_{L^{\infty}\left(B^{+}_{R}\right)}\left\lVert  D^{2}u \right\rVert_{\mathcal{L}^{^{(p,\theta)}}\left( \mathcal{U}_{2\sigma}\left(y'\right)\right)} \\ + C\left\lVert \nabla u \right\rVert_{\mathcal{L}^{^{(p,\theta)}}\left( \mathcal{U}_{2\sigma}\left(y'\right)\right)} + C\left\lVert  u \right\rVert_{\mathcal{L}^{^{(p,\theta)}}\left( \mathcal{U}_{2\sigma}\left(y'\right)\right)} \\+\left\lVert \widetilde{f} \right\rVert_{\mathcal{L}^{^{(p,\theta)}}\left( \mathcal{U}_{2\sigma}\left(y'\right)\right)}. 
			\end{multlined}\label{estimate for beta}
		\end{align}
		Now we write $\alpha=u-\beta.$  It is easy to see that $\beta \in W^{1,2}_{T}\left( \mathcal{U}_{R}; \varLambda^{k}\right)$ satisfies 
		\begin{multline}\label{zero BV on U}
			\int_{\mathcal{U}}\left<\bar{A}d\beta;d\psi\right> + \int_{\mathcal{U}}\left<d^{\ast}\left(\bar{B}\beta\right);d^{\ast}\left(\bar{B}\psi\right)\right>   + \int_{\mathcal{U}} \langle \widetilde{f} + \mathrm{P}u + \mathrm{R}\nabla u ; \psi \rangle 
			\\ + \int_{\mathcal{U}} \langle \mathrm{Q} u ; \nabla\psi \rangle   + \int_{\mathcal{U}} \langle \mathrm{S} \nabla u , \nabla\psi \rangle = 0 
		\end{multline}
		for every $\psi\in W^{1,2}_{T}\left( \mathcal{U}; \varLambda^{k}\right).$ Note that if we extend any $\psi \in W^{1,2}_{T}\left( \mathcal{U}_{R}; \varLambda^{k}\right)$ outside by zero, then $\psi \in W_{T, flat}^{1,2}(B_{R_{0}}^{+} ; \varLambda^{k}).$ Thus, from \eqref{flattened equation} and \eqref{zero BV on U}, we have 
		\begin{align}
			\int_{\mathcal{U}_{R}}\left\langle \bar{A}d\alpha;d\psi\right\rangle + \int_{\mathcal{U}_{R}}\left\langle d^{\ast}\left(\bar{B}\alpha\right);d^{\ast}\left(\bar{B}\psi\right)\right\rangle =0 
		\end{align}
		for all $\psi \in W_{T}^{1,2}(  \mathcal{U}_{R} ; \varLambda^{k})$ and $\alpha$ satisfies $\nu \wedge \alpha = 0$ on $\Gamma_{3R/4}.$ Note that by Lorentz regularity for \eqref{frozen system for beta} in smooth domains, we can easily show $\beta \in \mathsf{W}^{2}\mathcal{L}^{(p,\theta)}\left( \mathcal{U} ; \varLambda^{k}\right)$ and thus so is $\alpha.$  Now, for any $0 < \rho < \sigma/2,$ we have 
		\begin{align*}
			\left\lVert D^{2}u\right\rVert&_{\mathcal{L}^{\left(p, \theta \right)}\left( B^{+}_{2\rho}\left( y' \right) \right)} \\
			&\leq 	\left\lVert D^{2}\alpha\right\rVert_{\mathcal{L}^{\left(p, \theta \right)}\left( B^{+}_{2\rho}\left( y' \right) \right)} + \left\lVert D^{2}\beta\right\rVert_{\mathcal{L}^{\left(p, \theta \right)}\left( B^{+}_{2\rho}\left( y' \right) \right)}
			\\
			&\leq 	c\rho^{\left(\frac{n}{p}- \frac{n}{q}\right)}\left\lVert D^{2}\alpha\right\rVert_{\mathcal{L}^{\left(q, \theta \right)}\left( B^{+}_{2\rho}\left( y' \right) \right)} + \left\lVert D^{2}\beta\right\rVert_{\mathcal{L}^{\left(p, \theta \right)}\left( B^{+}_{2\rho}\left( y' \right) \right)}
			\\
			&\leq 	c\rho^{\left(\frac{n}{p}- \frac{n}{q}\right)}\left\lVert D^{2}\alpha\right\rVert_{\mathcal{L}^{\left(q, \theta \right)}\left( B^{+}_{\sigma}\left( y' \right) \right)} + \left\lVert D^{2}\beta\right\rVert_{\mathcal{L}^{\left(p, \theta \right)}\left( B^{+}_{2\sigma}\left( y' \right) \right)}
			\\
			&\leq 	c\left( \frac{\rho}{\sigma}\right)^{\left(\frac{n}{p}- \frac{n}{q}\right)}\left\lVert D^{2}\alpha\right\rVert_{\mathcal{L}^{\left(p, \theta \right)}\left( B^{+}_{3\sigma/2}\left( y' \right) \right)} + \left\lVert D^{2}\beta\right\rVert_{\mathcal{L}^{\left(p, \theta \right)}\left( B^{+}_{2\sigma}\left( y' \right) \right)}
			\\
			&\leq 	c\left( \frac{\rho}{\sigma}\right)^{\left(\frac{n}{p}- \frac{n}{q}\right)}\left\lVert D^{2}u\right\rVert_{\mathcal{L}^{\left(p, \theta \right)}\left( B^{+}_{2\sigma}\left( y' \right) \right)} + C\left\lVert D^{2}\beta\right\rVert_{\mathcal{L}^{\left(p, \theta \right)}\left( B^{+}_{2\sigma}\left( y' \right) \right)}. 
		\end{align*}
		Using this estimate and \eqref{estimate for beta}, we have 
		\begin{align*}
			\left\lVert D^{2}u\right\rVert_{\mathcal{L}^{\left(p, \theta \right)}\left( B^{+}_{2\rho}\left( y' \right) \right)} &\leq \begin{multlined}[t]
				C\left( \left( \frac{\rho}{\sigma}\right)^{\left(\frac{n}{p}- \frac{n}{q}\right)} + \left\lVert S \right\rVert_{L^{\infty}\left(B^{+}_{R}\right)} \right)\left\lVert  D^{2}u \right\rVert_{\mathcal{L}^{^{(p,\theta)}}\left( \mathcal{U}_{2\sigma}\left(y'\right)\right)} \\ + C\left\lVert \nabla u \right\rVert_{\mathcal{L}^{^{(p,\theta)}}\left( \mathcal{U}_{2\sigma}\left(y'\right)\right)} + C\left\lVert  u \right\rVert_{\mathcal{L}^{^{(p,\theta)}}\left( \mathcal{U}_{2\sigma}\left(y'\right)\right)} \\+\sigma^{\frac{\mu}{p}}\left\lVert \widetilde{f} \right\rVert_{\mathrm{L}_{\mu}^{(p,\theta)}\left( B^{+}_{R_{0}}\right)} . 
			\end{multlined}
		\end{align*}
		Now we claim that we have the estimate 
		\begin{align}\label{estimate for lower order terms}
			\left\lVert \nabla u \right\rVert_{\mathcal{L}^{^{(p,\theta)}}\left( \mathcal{U}_{2\sigma}\left(y'\right)\right)} + \left\lVert  u \right\rVert_{\mathcal{L}^{^{(p,\theta)}}\left( \mathcal{U}_{2\sigma}\left(y'\right)\right)} \leq C \sigma^{\frac{\tilde{\kappa}}{p}} \left\lVert u \right\rVert_{\mathsf{W}^{2}\mathrm{L}_{\kappa}^{(p,\theta)}\left( B^{+}_{R_{0}}\right)}. 
		\end{align}
		Note that since $\tilde{\kappa} \leq \mu < n,$ we can always choose $q$ large enough such that 
		$$ \frac{n}{p} - \frac{n}{q} > \frac{\tilde{\kappa}}{p} .$$ Thus, assuming the estimate \eqref{estimate for lower order terms}, we can use the standard iteration lemma ( Lemma 5.13 in \cite{giaquinta-martinazzi-regularity} ) to choose $R$ small enough, using \eqref{smallness of S}, such that $\left\lVert S \right\rVert_{L^{\infty}\left(B^{+}_{R}\right)}$ is smaller than the $\varepsilon_{0}$ given by the lemma. Then the lemma implies the estimate  
		\begin{align*}
			\rho^{-\frac{\tilde{\kappa}}{p}}\left\lVert D^{2}u\right\rVert&_{\mathcal{L}^{\left(p, \theta \right)}\left( B^{+}_{2\rho}\left( y' \right) \right)} \leq C \left( \left\lVert u \right\rVert_{\mathsf{W}^{l+2}\mathrm{L}_{\kappa}^{(p,\theta)}\left( B^{+}_{R_{0}}\right)}   +  \left\lVert \widetilde{f} \right\rVert_{\mathrm{L}_{\mu}^{(p,\theta)}\left( B^{+}_{R_{0}}\right)} \right). 
		\end{align*}\smallskip 
		
		Now only remains to prove the estimate \eqref{estimate for lower order terms}.  Note that by Morrey-Lorentz Sobolev embeddings, we have 
		\begin{align*}
			\mathsf{W}^{2}\mathrm{L}_{\kappa}^{(p,\theta)}\left( \mathcal{U}_{2\sigma}\left(y'\right)\right) \hookrightarrow \left\lbrace \begin{aligned}
				&\mathrm{L}_{\kappa}^{^{\left( \frac{\left(n-\kappa\right)p}{n - \kappa-2p},\frac{\left(n-\kappa\right)\theta}{n - \kappa-2p}\right)}}\left( \mathcal{U}_{2\sigma}\left(y'\right)\right) &&\text{ if } 2p < n-\kappa, \\
				&L^{r}\left( \mathcal{U}_{2\sigma}\left(y'\right)\right) &&\text{ if } 2p \geq n - \kappa, \\
				&\mathsf{W}^{1}\mathrm{L}_{\kappa}^{^{\left( \frac{\left(n-\kappa\right)p}{n - \kappa-2p},\frac{\left(n-\kappa\right)\theta}{n - \kappa-2p}\right)}}\left( \mathcal{U}_{2\sigma}\left(y'\right)\right) &&\text{ if } p < n-\kappa, \\
				&W^{1, r}\left( \mathcal{U}_{2\sigma}\left(y'\right)\right) &&\text{ if } p \geq n - \kappa, 
			\end{aligned}\right. 
		\end{align*}
		for any $1 \leq r < \infty.$ Now we estimate 
		\begin{align*}
			\left\lVert  u \right\rVert_{\mathcal{L}^{^{(p,\theta)}}\left( \mathcal{U}_{2\sigma}\left(y'\right)\right)} \leq \left\lbrace \begin{aligned}
				&C\sigma^{\frac{2p+\kappa}{p}}\left\lVert  u \right\rVert_{\mathrm{L}_{\kappa}^{^{\left( \frac{\left(n-\kappa\right)p}{n - \kappa-2p},\frac{\left(n-\kappa\right)\theta}{n - \kappa-2p}\right)}}\left( \mathcal{U}_{2\sigma}\left(y'\right)\right)} &&\text{ if } 2p < n-\kappa, \\
				&C\sigma^{n\left( \frac{1}{p} - \frac{1}{r}\right)}\left\lVert  u \right\rVert_{L^{r}\left( \mathcal{U}_{2\sigma}\left(y'\right)\right)} &&\text{ if } 2p \geq n - \kappa,
			\end{aligned}\right. 
		\end{align*}
		for any $1 \leq r < \infty.$ Similarly, we have 
		\begin{align*}
			\left\lVert  \nabla u \right\rVert_{\mathcal{L}^{^{(p,\theta)}}\left( \mathcal{U}_{2\sigma}\left(y'\right)\right)} \leq \left\lbrace \begin{aligned}
				&C\sigma^{\frac{p+\kappa}{p}}\left\lVert  u \right\rVert_{\mathrm{L}_{\kappa}^{^{\left( \frac{\left(n-\kappa\right)p}{n - \kappa-p},\frac{\left(n-\kappa\right)\theta}{n - \kappa-p}\right)}}\left( \mathcal{U}_{2\sigma}\left(y'\right)\right)} &&\text{ if } p < n-\kappa, \\
				&C\sigma^{n\left( \frac{1}{p} - \frac{1}{r}\right)}\left\lVert  u \right\rVert_{L^{r}\left( \mathcal{U}_{2\sigma}\left(y'\right)\right)} &&\text{ if } p \geq n - \kappa,
			\end{aligned}\right. 
		\end{align*} for any $1 \leq r < \infty.$ Now note that since by definition of $\tilde{\kappa},$ we have $\frac{\tilde{\kappa}}{p} \leq \frac{\mu}{p} < \frac{n}{p},$ we can choose $r$ large enough such that 
		\begin{align*}
			\frac{\tilde{\kappa}}{p} \leq	n\left( \frac{1}{p} - \frac{1}{r}\right)  < \frac{n}{p}.
		\end{align*} 
		This establishes \eqref{estimate for lower order terms} and completes the proof. \smallskip

		Now for $(ii),$ note that setting $\beta = B(x)\omega$ and $\psi = B(x)\phi,$ we see immediately that 
		it is enough to prove the regularity estimates for $\beta \in  W^{1,2}_{N}\left(\Omega; \Lambda^{k} \right)$ satisfying,  
		for all $\psi \in W^{1,2}_{N}\left(\Omega; \Lambda^{k} \right),$
		\begin{multline}\label{weakformW12N}
			\int_{\Omega}  \langle A (x)d \left( B^{-1}(x)\beta \right), d\left( B^{-1}(x)\psi \right) \rangle   + \int_{\Omega}\langle d^{\ast} \beta  , d^{\ast} \psi \rangle 
			\\ + \lambda  \int_{\Omega}\langle \left( B^{-1}(x)\right)^{T}\beta , \psi \rangle   
			+\int_{\Omega} \langle  f, \psi \rangle  = 0.   
		\end{multline}
		Now the estimate for this system follows in an analogous manner, this time using Remark \ref{flatteninglemmanormal}, Theorem \ref{BoundaryLorentzzeroBVPnormal} and $(ii)$ of Theorem \ref{boundary hessian decay}.   \end{proof}
	\subsection{Global estimates}		
	\begin{theorem}\label{global Hessian morrey-Lorentz estimate}
		Let $l \geq 0$ be an integer and let  $1 <p < \infty, 1 \leq \theta \leq \infty, 0 \leq \mu < n .$ Let $\Omega \subset \mathbb{R}^{n}$ be a open, bounded subset such that $\partial\Omega$ is of class $C^{l+2}.$ Let $A \in C^{l+1}\left(\overline{\Omega}; \operatorname{Hom}\left(\varLambda^{k+1}\right)\right)$  and $B \in C^{l+2}\left(\overline{\Omega}; \operatorname{Hom}\left(\varLambda^{k}\right)\right)$ satisfy the Legendre condition. Suppose  $f \in \mathsf{W}^{l}\mathrm{L}^{(p,\theta)}_{\mu}\left( \Omega; \varLambda^{k}\right) \cap W^{l,2}\left( \Omega; \varLambda^{k}\right)$ and $\lambda \in \mathbb{R}.$ Let $\omega \in W^{l+2,2}\left( \Omega; \varLambda^{k}\right) \cap \mathsf{W}^{l+2}\mathcal{L}^{(p,\theta)}\left(\Omega; \varLambda^{k}\right)$ and assume one of the following holds.\smallskip 
		
		\begin{enumerate}[(i)]
			\item Let $\omega \in W_{T}^{1,2}(\Omega, \varLambda^{k})$ and for all $ \phi \in W_{T}^{1,2}(\Omega; \varLambda^{k}),$ we have,  
			\begin{multline}\label{bilinear form global estimate}
				\int_{\Omega}  \langle A (x)d\omega , d\phi \rangle   + \int_{\Omega}\langle d^{\ast} \left( B(x)\omega \right) , d^{\ast} \left( B(x)\phi \right) \rangle 
				+ \lambda  \int_{\Omega}\langle  B(x)\omega , \phi \rangle 	\\ +\int_{\Omega} \langle  f, \phi \rangle  = 0.    
			\end{multline}
			
			\item Let $\omega \in W_{B,N}^{1,2}(\Omega, \varLambda^{k})$ and for all $ \phi \in W_{B,N}^{1,2}(\Omega; \varLambda^{k}),$ we have, 
			\begin{multline}\label{bilinear form global estimate normal}
				\int_{\Omega}  \langle A (x)d\omega , d\phi \rangle   + \int_{\Omega}\langle d^{\ast} \left( B(x)\omega \right) , d^{\ast} \left( B(x)\phi \right) \rangle 
				+ \lambda  \int_{\Omega}\langle  B(x)\omega , \phi \rangle 	\\	+\int_{\Omega} \langle  f, \phi \rangle  = 0.    
			\end{multline}
		\end{enumerate}
		\noindent  Then $\omega \in \mathsf{W}^{l+2}\mathrm{L}^{(p,\theta)}_{\mu}\left( \Omega; \varLambda^{k}\right)$ and there exists a constant $$C = C \left( l, n, k, N, p, \theta, \mu, \Omega, \left\lVert A \right\rVert_{C^{l+1}}, \left\lVert B \right\rVert_{C^{l+2}}, \gamma, \lambda  \right) >0$$ such that 
		\begin{align*}
			\left\lVert u \right\rVert_{\mathsf{W}^{l+2}\mathrm{L}^{(p,\theta)}_{\mu}\left( \Omega \right)} \leq C \left( \left\lVert u \right\rVert_{\mathsf{W}^{l+2}\mathcal{L}^{(p,\theta)}\left( \Omega\right)} + \left\lVert f \right\rVert_{\mathsf{W}^{l}\mathrm{L}^{(p,\theta)}_{\mu}\left( \Omega\right)} \right). 
		\end{align*}
	\end{theorem}
	\begin{proof}
		We only show $(i).$ We prove the case $l=0,$ as the argument can be iterated. By standard localization, covering and gluing argument, the result can be deduced from the local interior and boundary estimates. Since the interior estimate follows the same way and are much easier, the boundary estimates of Lemma \ref{boundary estimate lemma} implies 
		that if $\omega \in \mathsf{W}^{l}\mathrm{L}^{(p,\theta)}_{\kappa}\left( \Omega\right)$ for any $0 \leq \kappa < \mu ,$ then we have the estimate 
		\begin{align*}
			\left\lVert \omega \right\rVert_{\mathsf{W}^{l+2}\mathrm{L}^{(p,\theta)}_{\tilde{\kappa}}\left( \Omega \right)} \leq C \left( \left\lVert \omega \right\rVert_{\mathsf{W}^{l+2}\mathrm{L}_{\kappa}^{(p,\theta)}\left( \Omega\right)} + \left\lVert f \right\rVert_{\mathsf{W}^{l}\mathrm{L}^{(p,\theta)}_{\mu}\left( \Omega\right)} \right), 
		\end{align*}  where $\tilde{\kappa} = \min \left\lbrace \kappa + p, \mu \right\rbrace.$ Since we have assumed $\omega \in \mathsf{W}^{l+2}\mathcal{L}^{(p,\theta)}\left( \Omega\right),$ we can start from $\kappa =0$ and bootstrap. Now since $\mu < n$ and $p >1,$ there exists a natural number $ 1 \leq m \leq n$ such that $$ \left(m-1\right)p \leq \mu < mp. $$ Thus, we would establish our desired estimate in at most $m$ steps.
	\end{proof}
	\subsection{Approximation}
	We would need an approximation result, which is somewhat non-standard, as our spaces are in general not separable.  
	\begin{lemma}\label{approximation lemma}[Approximation Lemma]
		Suppose for some $\lambda \in \mathbb{R},$ there exists unique solution $\omega \in \mathsf{W}^{l+2}\mathrm{L}^{(p,\theta)}_{\mu}\left( \Omega; \varLambda^{k}\right) \cap W^{l+2,2}\left( \Omega; \varLambda^{k}\right)$ to \eqref{bvp hodge elliptic system general full regularity} ( respectively, \eqref{bvp hodge elliptic system general full regularity normal} ) for any given $f \in \mathsf{W}^{l}\mathrm{L}^{(p,\theta)}_{\mu}\left( \Omega; \varLambda^{k}\right) \cap W^{l,2}\left( \Omega; \varLambda^{k}\right)$  and any given $\omega_{0} \in \mathsf{W}^{l+2}\mathrm{L}^{(p,\theta)}_{\mu}\left( \Omega; \varLambda^{k}\right) \cap W^{l+2,2}\left( \Omega; \varLambda^{k}\right)$ which satisfies the estimate 
		\begin{align}\label{apriori estimate}
			\left\lVert \omega \right\rVert_{\mathsf{W}^{l+2}\mathrm{L}^{(p,\theta)}_{\mu}\left( \Omega \right)} \leq C \left( \left\lVert f \right\rVert_{\mathsf{W}^{l}\mathrm{L}^{(p,\theta)}_{\mu}\left( \Omega\right)}  + \left\lVert \omega_{0} \right\rVert_{\mathsf{W}^{l+2}\mathrm{L}^{(p,\theta)}_{\mu}\left( \Omega\right)}\right).	
		\end{align}
		Then for any $f \in \mathsf{W}^{l}\mathrm{L}^{(p,\theta)}_{\mu}\left( \Omega; \varLambda^{k}\right) $ and $\omega_{0} \in \mathsf{W}^{l+2}\mathrm{L}^{(p,\theta)}_{\mu}\left( \Omega; \varLambda^{k}\right),$ there exists $\omega \in \mathsf{W}^{l+2}\mathrm{L}^{(p,\theta)}_{\mu}\left( \Omega; \varLambda^{k}\right)$  which satisfies the estimate \eqref{apriori estimate} and is the unique solution to \eqref{bvp hodge elliptic system general full regularity} ( respectively, \eqref{bvp hodge elliptic system general full regularity normal} ). 
	\end{lemma}
	\begin{proof}
		We show the case for \eqref{bvp hodge elliptic system general full regularity} and prove only for $l=0.$ By considering the system for $\omega - \omega_{0},$ we can assume $\omega_{0} = 0.$ Now for any $f \in \mathrm{L}^{(p,\theta)}_{\mu}\left( \Omega; \varLambda^{k}\right),$ by extension and mollification, it is easy to check using Young's inequality for convolutions in Lorentz spaces ( proved by O'Neil in \cite{ONeil_ConvolutionLorentzSpaces} ) and Jensen's inequality that we can find a sequence $\left\lbrace f_{s}\right\rbrace_{s \in \mathbb{N}} \subset C^{\infty}\left( \overline{\Omega}; \varLambda^{k}\right)$ such that we have 
		\begin{align}\label{approximation by smooth functions}
			\left\lbrace \begin{aligned}
				\limsup\limits_{s \in \mathbb{N}} \left\lVert f_{s} \right\rVert_{\mathrm{L}^{(p,\theta)}_{\mu}\left( \Omega\right)} &\leq C \left\lVert f \right\rVert_{\mathrm{L}^{(p,\theta)}_{\mu}\left( \Omega\right)} \\
				f_{s} &\rightarrow f \qquad \text{ strongly in } L^{q}\left( \Omega; \varLambda^{k}\right),
			\end{aligned}\right. 
		\end{align}
		where $q>1$ is any exponent such that $\mathcal{L}^{\left(p, \theta\right)} \subset L^{q}.$ Note that unless $\mu = 0$ and $\theta \neq \infty,$ \eqref{approximation by smooth functions} can not be improved to strong convergence in $\mathrm{L}^{(p,\theta)}_{\mu}$ norms, as the corresponding spaces are not separable and smooth functions are not dense. However, this would be good enough. Indeed, by using the hypothesis, there exists a sequence $\left\lbrace \omega_{s}\right\rbrace_{s \in \mathbb{N}}$ such that for each $s \in \mathbb{N},$ $\omega_{s}$ is the unique solution to 
		\begin{align}\label{bvp hodge approximation lemma for each s}
			\left\lbrace \begin{gathered}
				d^{\ast} ( A (x) d\omega_{s} )  + \left( B(x) \right)^{T}d d^{\ast}\left( B(x) \omega_{s} \right)   =  \lambda B(x)\omega_{s} + f_{s}  \text{ in } \Omega, \\
				\nu\wedge \omega_{s} = 0 \text{  on } \partial\Omega, \\
				\nu\wedge d^{\ast} \left( B(x)\omega_{s} \right) = 0 \text{ on } \partial\Omega, 
			\end{gathered} 
			\right. 
		\end{align}
		and satisfies the estimate 
		\begin{align}\label{apriori estimate for each s}
			\left\lVert \omega_{s} \right\rVert_{\mathsf{W}^{2}\mathrm{L}^{(p,\theta)}_{\mu}\left( \Omega \right)} \leq C  \left\lVert f_{s} \right\rVert_{\mathrm{L}^{(p,\theta)}_{\mu}\left( \Omega\right)} \leq C \left\lVert f \right\rVert_{\mathrm{L}^{(p,\theta)}_{\mu}\left( \Omega\right)} .	
		\end{align}
		This implies, in particular, that $\left\lbrace \omega_{s}\right\rbrace_{s \in \mathbb{N}}$ is uniformly bounded in $\mathsf{W}^{2}\mathcal{L}^{\left(p, \theta\right)}.$ Since these spaces are reflexive when $1 < \theta < \infty,$ these immediately imply that up to the extraction of a subsequence that we do not relabel, there exist $\omega \in \mathsf{W}^{2}\mathcal{L}^{\left(p, \theta\right)}$ such that 
		\begin{align*}
			\omega_{s} \rightharpoonup \omega \qquad \text{ weakly in }  \mathsf{W}^{2}\mathcal{L}^{\left(p, \theta\right)}. 
		\end{align*}
		Now it is easy to show that this and the uniform bound in $ \mathsf{W}^{2}\mathrm{L}^{(p,\theta)}_{\mu}$ implies 
		\begin{align*}
			\left\lVert \omega \right\rVert_{\mathsf{W}^{2}\mathrm{L}^{(p,\theta)}_{\mu}\left( \Omega \right)} \leq \liminf\limits_{s\rightarrow \infty}\left\lVert \omega_{s} \right\rVert_{\mathsf{W}^{2}\mathrm{L}^{(p,\theta)}_{\mu}\left( \Omega \right)} \leq C \left\lVert f \right\rVert_{\mathrm{L}^{(p,\theta)}_{\mu}\left( \Omega\right)}. 
		\end{align*}
		Since the weak convergence in $\mathsf{W}^{2}\mathcal{L}^{\left(p, \theta\right)}$ implies weak convergence in $W^{2,q}$ for any $1 < q<p,$ which allows us to pass to the limit in \eqref{bvp hodge approximation lemma for each s} to conclude $\omega$ solves \eqref{bvp hodge elliptic system general full regularity}. Replacing weak convergence by weak star convergence, this argument works also when $\theta = \infty.$ In this case, the uniform bound implies, by virtue of separability of $\mathcal{L}^{\left(p, 1\right)},$ that we have  
		\begin{align*}
			\omega_{s} \stackrel{\ast}{\rightharpoonup} \omega \qquad \text{ weakly $\ast$ in }  \mathsf{W}^{2}\mathcal{L}^{\left(p, \infty\right)}. 
		\end{align*}
		The rest follows exactly as before. Indeed, by the weak star convergence, we have 
		\begin{align*}
			\left\lVert \omega \right\rVert_{\mathsf{W}^{2}\mathrm{L}^{(p,\infty)}_{\mu}\left( \Omega \right)} \leq \liminf\limits_{s\rightarrow \infty}\left\lVert \omega_{s} \right\rVert_{\mathsf{W}^{2}\mathrm{L}^{(p,\infty)}_{\mu}\left( \Omega \right)} \leq C \left\lVert f \right\rVert_{\mathrm{L}^{(p,\infty)}_{\mu}\left( \Omega\right)}. 
		\end{align*}
		The fact that $\omega$ solves \eqref{bvp hodge elliptic system general full regularity} can be checked as before. Thus, it remains to settle the case $\theta =1,$ which is trickier. The uniform bound in $\mathsf{W}^{2}\mathcal{L}^{\left(p, 1\right)}$ implies a uniform bound in $\mathsf{W}^{2}\mathcal{L}^{\left(p, \tilde{\theta}\right)}$ for any $1 < \tilde{\theta} < \infty.$ Thus, up to the extraction of a subsequence that we do not relabel, there exist $\omega \in \mathsf{W}^{2}\mathcal{L}^{\left(p, \theta\right)}$ such that 
		\begin{align*}
			\omega_{s} \rightharpoonup \omega \qquad \text{ weakly in }  \mathsf{W}^{2}\mathcal{L}^{\left(p, \tilde{\theta}\right)} \qquad \text{ for every } 1 < \tilde{\theta} < \infty.
		\end{align*}
		But we can argue using Fatou's lemma as in the proof of Theorem 7.1 in \cite{Costea_SobolevLorentzCapacity} to conclude that the above weak convergence and the uniform bound in  $\mathsf{W}^{2}\mathcal{L}^{\left(p, 1\right)}$ implies $\omega \in \mathsf{W}^{2}\mathcal{L}^{\left(p, 1\right)}$ and we have 
		\begin{align*}
			\left\lVert \omega \right\rVert_{\mathsf{W}^{2}\mathcal{L}^{(p,1)}\left( \Omega \right)} \leq \liminf\limits_{s\rightarrow \infty}\left\lVert \omega_{s} \right\rVert_{\mathsf{W}^{2}\mathcal{L}^{(p,1)}\left( \Omega \right)}. 
		\end{align*}
		But note that this last argument holds as well when we restrict everything to any open subset of $\Omega$ as well. So using this for $\Omega_{(\rho , x_{0})}$ for $x_{0} \in \Omega$ and $\rho >0,$ we get 
		\begin{align*}
			\rho^{- \frac{\mu}{p}}\left\lVert \omega \right\rVert_{\mathsf{W}^{2}\mathcal{L}^{(p,1)}\left( \Omega_{(\rho , x_{0})} \right)} \leq 	\rho^{- \frac{\mu}{p}}\liminf\limits_{s\rightarrow \infty}\left\lVert \omega_{s} \right\rVert_{\mathsf{W}^{2}\mathcal{L}^{(p,1)}\left( \Omega_{(\rho , x_{0})} \right)}. 
		\end{align*}
		Taking supremum over $\rho >0$ and $x_{0} \in \Omega,$ we deduce 
		\begin{align*}
			\left\lVert \omega \right\rVert_{\mathsf{W}^{2}\mathrm{L}^{(p,1)}_{\mu}\left( \Omega \right)} \leq \liminf\limits_{s\rightarrow \infty}\left\lVert \omega_{s} \right\rVert_{\mathsf{W}^{2}\mathrm{L}^{(p,1)}_{\mu}\left( \Omega \right)} \leq C \left\lVert f \right\rVert_{\mathrm{L}^{(p,1)}_{\mu}\left( \Omega\right)}. 
		\end{align*}
		Once again, it is easy to check $\omega$ solves \eqref{bvp hodge elliptic system general full regularity}. Uniqueness of $\omega$, in all cases, follow from the fact that if $\omega_{1}$ and $\omega_{2}$ are two solutions, then $\omega_{1} - \omega_{2}$ satisfies \eqref{bvp hodge elliptic system general full regularity} with $f =0$ and $\omega_{0}=0.$ Then the uniqueness assumption in the hypothesis implies $\omega_{1}- \omega_{2}=0.$ This completes the proof. 
	\end{proof}
	\section{Main results}\label{main results}
	Throughout this entire section, we would assume, without specific mention, that $n \geq 2,$ $N \geq 1,$ $1 \leq k \leq n-1,$  $l \geq 0$ are integers and $\Omega \subset \mathbb{R}^n$ be an open, bounded $C^{l+2}$ set. The exponents $p, \theta, \mu $ satisfy $1 < p < \infty$, $1\leq\theta \leq \infty$, $0\leq\mu<n.$
	
	\subsection{Morrey-Lorentz estimate for Hodge systems}
	\begin{theorem}\label{generalHodgesystemtheorem}
		Let $A \in C^{l+1}\left(\overline{\Omega}; \operatorname{Hom}\left(\varLambda^{k+1}\right)\right)$ and $B \in C^{l+2}\left(\overline{\Omega}; \operatorname{Hom}\left(\varLambda^{k}\right)\right)$ satisfy the Legendre condition.
		Then the following holds.
		\begin{enumerate}[(i)]
			\item There exists an at most countable set 
			$\sigma \subset ( -\infty, 0]$, with no limit points except possibly $- \infty,$ such that the following boundary value problem,
			\begin{equation}\label{eigenvalue problem hodge system general full regularity}
				\left\lbrace \begin{gathered}
					d^{\ast} ( A (x) d\alpha ) + \left( B(x) \right)^{T}d d^{\ast}\left( B(x) \alpha \right)  = \sigma_{i} B(x) \alpha  \text{ in } \Omega, \\
					\nu\wedge \alpha = 0 \text{  on } \partial\Omega, \\
					\nu\wedge d^{\ast} \left( B(x) \alpha \right) = 0 \text{ on } \partial\Omega.
				\end{gathered} 
				\right. \tag{$\mathbb{EP}_{T}$}
			\end{equation}
			has non-trivial solutions $\alpha$ if and only if  $\sigma_{i} \in \sigma.$ For any $1 < \tilde{p}< \infty,$ $1 \leq \tilde{\theta} \leq \infty$ and  $0 \leq \tilde{\mu} < n,$ all such solutions $\alpha \in \mathsf{W}^{l+2}\mathrm{L}^{(\tilde{p},\tilde{\theta})}_{\tilde{\mu}}\left( \Omega; \varLambda^{k}\right).$ Also, for any $\sigma_{i} \in \sigma,$ the space of solutions to \eqref{eigenvalue problem hodge system general full regularity}, denoted $\mathcal{E}_{i, T}$ is a finite-dimensional subspace of $\mathrm{L}^{(p,\theta)}_{\mu}\left( \Omega; \varLambda^{k}\right)$   and  $\operatorname{dim}\mathcal{E}_{i, T} = \operatorname{dim}\mathcal{E}_{i, T}^{\ast},$ where $\mathcal{E}_{i, T}^{\ast}$ denotes the space of solutions of 
			\begin{equation}\label{bvp hodge elliptic system general full regularity adjoint}
				\left\lbrace \begin{gathered}
					d^{\ast} ( \left( A(x) \right)^{\intercal} d\psi )  + \left( B(x) \right)^{\intercal}d d^{\ast}\left( B(x) \psi \right)   =  \sigma_{i} \left( B(x) \right)^{\intercal}\psi  \text{ in } \Omega, \\
					\nu\wedge \psi = 0 \text{  on } \partial\Omega, \\
					\nu\wedge d^{\ast} \left( B(x)\psi \right) = 0 \text{ on } \partial\Omega. 
				\end{gathered} 
				\right. \tag{$\mathbb{EP}_{T}^{\ast}$}
			\end{equation}
			\item If $\lambda \notin \sigma$, then 
			for any $f \in \mathsf{W}^{l}\mathrm{L}^{(p,\theta)}_{\mu}\left( \Omega; \varLambda^{k}\right)$, and any $\omega_{0} \in \mathsf{W}^{l+2}\mathrm{L}^{(p,\theta)}_{\mu}\left( \Omega; \varLambda^{k}\right) $ there exists a unique 
			solution $\omega \in \mathsf{W}^{l+2}\mathrm{L}^{(p,\theta)}_{\mu}\left( \Omega; \varLambda^{k}\right)$ to the following boundary value problem:
			\begin{equation}\label{bvp hodge elliptic system general full regularity}
				\left\lbrace \begin{gathered}
					d^{\ast} ( A (x) d\omega )  + \left( B(x) \right)^{T}d d^{\ast}\left( B(x) \omega \right)   =  \lambda B(x)\omega + f  \text{ in } \Omega, \\
					\nu\wedge \omega = \nu\wedge\omega_{0} \text{  on } \partial\Omega, \\
					\nu\wedge d^{\ast} \left( B(x)\omega \right) = \nu\wedge d^{\ast} \left( B(x) \omega_{0} \right) \text{ on } \partial\Omega, 
				\end{gathered} 
				\right. \tag{$\mathbb{P}_{T}$}
			\end{equation}
			which satisfies the estimate
			\begin{align}\label{full Hodge system estimate}
				\left\lVert \omega \right\rVert_{\mathsf{W}^{l+2}\mathrm{L}^{(p,\theta)}_{\mu}\left( \Omega \right)} \leq C \left( \left\lVert \omega \right\rVert_{\mathsf{W}^{l+2}\mathcal{L}^{(p,\theta)}\left( \Omega\right)} + \left\lVert f \right\rVert_{\mathsf{W}^{l}\mathrm{L}^{(p,\theta)}_{\mu}\left( \Omega\right)} + \left\lVert \omega_{0} \right\rVert_{\mathsf{W}^{l+2}\mathrm{L}^{(p,\theta)}_{\mu}\left( \Omega\right)}\right). 
			\end{align}
			\item If $\lambda = \sigma_{i}$ for some $i \in \mathbb{N},$ then for any $\omega_{0} \in \mathsf{W}^{l+2}\mathrm{L}^{(p,\theta)}_{\mu}\left( \Omega; \varLambda^{k}\right) $ and any $f \in \mathsf{W}^{l}\mathrm{L}^{(p,\theta)}_{\mu}\left( \Omega; \varLambda^{k}\right)$ satisfying 
			\begin{align*}
				\int_{\Omega} \left\langle f , \psi \right\rangle = 0  \qquad \text{ for all } \psi \in \mathcal{E}_{i, T}^{\ast},
			\end{align*}
			there exists a  unique solution $\omega \in \mathsf{W}^{l+2}\mathrm{L}^{(p,\theta)}_{\mu}\left( \Omega; \varLambda^{k}\right)/\mathcal{E}_{i, T}$ to \eqref{bvp hodge elliptic system general full regularity} satisfying estimate \eqref{full Hodge system estimate}
		\end{enumerate}
	\end{theorem}
	\begin{remark}\label{remarkhodgegeneraltan}
		If $\mathcal{H}_{T}\left( \Omega; \varLambda^{k} \right) \neq \lbrace 0 \rbrace,$ then it can be proved that $\alpha $ is a nontrivial solution for \eqref{bvp hodge elliptic system general full regularity adjoint} with 
		$\sigma_{i} = 0$ if and only if $\alpha = d\beta + h,$ where 
		where $\beta$ is a solution of 
		\begin{align*}
			\left\lbrace \begin{aligned} d^{\ast} ( B d \beta ) &=  - d^{\ast} ( Bh ) &&\text{ in } \Omega, \\
				d^{\ast} \beta &= 0 &&\text{ in } \Omega, \\
				\nu \wedge \beta &= 0 &&\text{ on  } \partial \Omega,\end{aligned}\right. \end{align*} for some nontrivial 
		$h \in \mathcal{H}_{T}\left( \Omega; \varLambda^{k} \right).$  Note also that if $B$ is a constant multiple of the identity matrix, then $\beta \in \mathcal{H}_{T}\left( \Omega; \varLambda^{k-1}\right)$ and thus $d\beta = 0.$ Consequently, \eqref{bvp hodge elliptic system general full regularity} with $\lambda = 0$ can be solved for any $f$ satisfying $ f \in \left( \mathcal{H}_{T}\left( \Omega; \varLambda^{k} \right)\right)^{\perp},$ if $B \equiv c\mathbb{I}$ for some constant $c>0.$ For a general $B,$ an additional condition $d^{\ast} f = 0$ in $\Omega$ is needed. If $\mathcal{H}_{T}\left( \Omega; \varLambda^{k} \right) = \lbrace 0 \rbrace,$ then \eqref{bvp hodge elliptic system general full regularity} with $\lambda = 0$ can be solved for any $f,$ no extra condition on $f$ is needed.    
	\end{remark}
	\begin{proof}
		Note that if $\mathcal{L}^{\left(p, \theta\right)} \subset L^{2},$ then standard Lax-Milgram theorem argument proves the existence of a $W^{1,2}$ weak solution ( see \cite{Sil_linearregularity} ) for $\lambda >0.$  Then Theorem \ref{BoundaryLorentzzeroBVP} applied to $\omega - \omega_{0}$ establishes the $\mathsf{W}^{l+2}\mathcal{L}^{\left(p, \theta\right)}$ estimates. Thus, we can then use Theorem \ref{global Hessian morrey-Lorentz estimate} to establish Morrey-Lorentz estimates. The approximation argument in Theorem \ref{approximation lemma} extends these estimates to the cases when $\mathcal{L}^{\left(p, \theta\right)} \not\subset L^{2}$ as well. The rest is standard Riesz-Fredholm theory, by virtue of the compact embedding $\mathsf{W}^{1}\mathcal{L}^{\left(p, \theta\right)} \hookrightarrow \mathcal{L}^{\left(p, \theta\right)}.$ This finishes the proof. 
	\end{proof}
	In an analogous manner, we have the following. 
	\begin{theorem}\label{generalHodgesystemtheoremnormal}
		Let $A \in C^{l+1}\left(\overline{\Omega}; \operatorname{Hom}\left(\varLambda^{k+1}\right)\right)$, and $B \in C^{l+2}\left(\overline{\Omega}; \operatorname{Hom}\left(\varLambda^{k}\right)\right)$, both satisfy the Legendre condition.
		Then the following holds.
		\begin{enumerate}
			\item There exists an at most countable set 
			$\sigma \subset ( -\infty, 0]$, with no limit points except possibly $- \infty,$ such that the following boundary value problem,
			\begin{equation}\label{eigenvalue problem hodge system general full regularity normal}
				\left\lbrace \begin{gathered}
					d^{\ast} ( A (x) d\alpha ) + \left( B(x) \right)^{T}d d^{\ast}\left( B(x) \alpha \right)  = \sigma_{i} B(x) \alpha  \text{ in } \Omega, \\
					\nu\lrcorner \left( B(x)\alpha \right) = 0 \text{  on } \partial\Omega, \\
					\nu\lrcorner \left( A(x) d\alpha \right) = 0 \text{ on } \partial\Omega.
				\end{gathered} 
				\right. \tag{$\mathbb{EP}_{N}$}
			\end{equation}
			has non-trivial solutions $\alpha \in \mathsf{W}^{l+2}\mathrm{L}^{(p,\theta)}_{\mu}\left( \Omega; \varLambda^{k}\right)$ if and only if  $\sigma_{i} \in \sigma.$ For any $1 < \tilde{p}< \infty,$ $1 \leq \tilde{\theta} \leq \infty$ and  $0 \leq \tilde{\mu} < n,$ all such solutions $\alpha \in \mathsf{W}^{l+2}\mathrm{L}^{(\tilde{p},\tilde{\theta})}_{\tilde{\mu}}\left( \Omega; \varLambda^{k}\right).$ Also, for any $\sigma_{i} \in \sigma,$ the space of solutions to \eqref{eigenvalue problem hodge system general full regularity normal}, denoted $\mathcal{E}_{i, N}$ is a finite-dimensional subspace of $\mathrm{L}^{(p,\theta)}_{\mu}\left( \Omega; \varLambda^{k}\right)$   and  $\operatorname{dim}\mathcal{E}_{i, N} = \operatorname{dim}\mathcal{E}_{i, N}^{\ast},$ where $\mathcal{E}_{i, N}^{\ast}$ denotes the space of solutions of 
			\begin{equation}\label{bvp hodge elliptic system general normal full regularity adjoint}
				\left\lbrace \begin{gathered}
					d^{\ast} ( \left( A(x) \right)^{\intercal} d\psi )  + \left( B(x) \right)^{\intercal}d d^{\ast}\left( B(x) \psi \right)   =  \sigma_{i} \left( B(x) \right)^{\intercal}\psi  \text{ in } \Omega, \\
					\nu\lrcorner \left( B(x) \psi \right) = 0 \text{  on } \partial\Omega, \\
					\nu\lrcorner \left( \left( A(x) \right)^{\intercal} d\psi \right)  = 0 \text{ on } \partial\Omega. 
				\end{gathered} 
				\right. \tag{$\mathbb{EP}_{N}^{\ast}$}
			\end{equation} 
			\item If $\lambda \notin \sigma$, then 
			for any $f \in \mathsf{W}^{l}\mathrm{L}^{(p,\theta)}_{\mu}\left( \Omega; \varLambda^{k}\right)$, and any $\omega_{0} \in \mathsf{W}^{l+2}\mathrm{L}^{(p,\theta)}_{\mu}\left( \Omega; \varLambda^{k}\right),$ there exists a unique 
			solution $\omega \in \mathsf{W}^{l+2}\mathrm{L}^{(p,\theta)}_{\mu}\left( \Omega; \varLambda^{k}\right),$ to the following boundary value problem:
			\begin{equation}\label{bvp hodge elliptic system general full regularity normal}
				\left\lbrace \begin{gathered}
					d^{\ast} ( A (x) d\omega )  + \left( B(x) \right)^{T}d d^{\ast}\left( B(x) \omega \right)   =  \lambda B(x)\omega + f  \text{ in } \Omega, \\
					\nu\lrcorner \left( B(x)\omega \right) = \nu\lrcorner\left( B(x)\omega_{0} \right) \text{  on } \partial\Omega. \\
					\nu\lrcorner  \left( A(x)d\omega \right) = \nu\lrcorner \left( A(x) d\omega_{0} \right) \text{ on } \partial\Omega, 
				\end{gathered} 
				\right. \tag{$\mathbb{P}_{N}$}
			\end{equation}
			which satisfies the estimate
			\begin{align}\label{full Hodge system estimate normal}
				\left\lVert \omega \right\rVert_{\mathsf{W}^{l+2}\mathrm{L}^{(p,\theta)}_{\mu}\left( \Omega \right)} \leq C \left( \left\lVert \omega \right\rVert_{\mathsf{W}^{l+2}\mathcal{L}^{(p,\theta)}\left( \Omega\right)} + \left\lVert f \right\rVert_{\mathsf{W}^{l}\mathrm{L}^{(p,\theta)}_{\mu}\left( \Omega\right)} + \left\lVert \omega_{0} \right\rVert_{\mathsf{W}^{l+2}\mathrm{L}^{(p,\theta)}_{\mu}\left( \Omega\right)}\right). 
			\end{align}
			\item If $\lambda = \sigma_{i}$ for some $i \in \mathbb{N},$ then for any $\omega_{0} \in \mathsf{W}^{l+2}\mathrm{L}^{(p,\theta)}_{\mu}\left( \Omega; \varLambda^{k}\right) $ and any $f \in \mathsf{W}^{l}\mathrm{L}^{(p,\theta)}_{\mu}\left( \Omega; \varLambda^{k}\right)$ satisfying 
			\begin{align*}
				\int_{\Omega} \left\langle f , \psi \right\rangle = 0  \qquad \text{ for all } \psi \in \mathcal{E}_{i, N}^{\ast},
			\end{align*}
			there exists a  unique solution $\omega \in \mathsf{W}^{l+2}\mathrm{L}^{(p,\theta)}_{\mu}\left( \Omega; \varLambda^{k}\right)/\mathcal{E}_{i, N}$ to \eqref{bvp hodge elliptic system general full regularity normal} satisfying estimate \eqref{full Hodge system estimate normal}
		\end{enumerate}
	\end{theorem}
	\begin{remark}\label{remarkhodgegeneralnormal}
		Analogously to Remark \ref{remarkhodgegeneraltan}, if $\mathcal{H}_{N}\left( \Omega; \varLambda^{k} \right) \neq \lbrace 0 \rbrace,$ then it can be proved that $\psi $ is a nontrivial solution for \eqref{bvp hodge elliptic system general normal full regularity adjoint} with $\sigma_{i} = 0$ if and only if $B^{-1}\psi = d^{\ast}\beta + h,$ where 
		where $\beta$ is a solution of 
		\begin{align*}
			\left\lbrace \begin{aligned} d ( B^{-1} d^{\ast} \beta ) &=  - d ( B^{-1}h ) &&\text{ in } \Omega, \\
				d\beta &= 0 &&\text{ in } \Omega, \\
				\nu \lrcorner \beta &= 0 &&\text{ on  } \partial \Omega,\end{aligned}\right. \end{align*} for some nontrivial 
		$h \in \mathcal{H}_{N}\left( \Omega; \varLambda^{k} \right).$  Thus if $B$ is a constant multiple of the identity matrix, then $\beta \in \mathcal{H}_{N}\left( \Omega; \varLambda^{k-1}\right)$ and thus $d^{\ast}\beta = 0.$ Consequently, \eqref{bvp hodge elliptic system general full regularity normal} with $\lambda = 0$ can be solved for any $f$ satisfying $ f \in \left( \mathcal{H}_{N}\left( \Omega; \varLambda^{k} \right)\right)^{\perp},$ if $B \equiv c\mathbb{I}$ for some constant $c>0.$ For a general $B,$ an additional condition $d \left( \left[ B^{-1}\right]^{\intercal}f\right) = 0$ in $\Omega$ is needed. If $\mathcal{H}_{N}\left( \Omega; \varLambda^{k} \right) = \lbrace 0 \rbrace,$ then \eqref{bvp hodge elliptic system general full regularity} with $\lambda = 0$ can be solved for any $f,$ no extra condition on $f$ is needed.    
	\end{remark}
	\subsection{Hodge decomposition in Morrey-Lorentz spaces}
	\begin{theorem}[Hodge decomposition]\label{Hodge decomposition}
		Let $A \in C^{l+1}\left(\overline{\Omega}; \operatorname{Hom}\left(\varLambda^{k+1}\right)\right)$ satisfy the Legendre condition and let $f \in \mathsf{W}^{l}\mathrm{L}^{(p,\theta)}_{\mu}\left( \Omega; \varLambda^{k}\right).$ Then the following holds. 
		\begin{itemize}
			\item[(i)] There exist $\alpha \in \mathsf{W}^{l+1}\mathrm{L}^{(p,\theta)}_{\mu}\left( \Omega; \varLambda^{k-1}\right)$ and $\beta \in \mathsf{W}^{l+1}\mathrm{L}^{(p,\theta)}_{\mu}\left( \Omega; \varLambda^{k+1}\right)$ and $h \in \mathcal{H}_{T}\left(\Omega; \varLambda^{k}\right)$ such that 
			\begin{gather*}
				f = d\alpha  + d^{\ast}\left( A\left(x\right)\beta\right) + h \qquad \text{ in } \Omega, \\
				d^{\ast}\alpha = d\beta = 0 \quad \text{ in } \Omega, \qquad 
				\nu\wedge \alpha  = \nu \wedge \beta = 0 \quad \text{ on } \partial\Omega. 
			\end{gather*}
			Moreover, we have the estimate 
			\begin{align*}
				\left\lVert \alpha \right\rVert_{\mathsf{W}^{l+1}\mathrm{L}^{(p,\theta)}_{\mu}} + 	\left\lVert \beta \right\rVert_{\mathsf{W}^{l+2}\mathrm{L}^{(p,\theta)}_{\mu}} + 	\left\lVert h \right\rVert_{\mathsf{W}^{l}\mathrm{L}^{(p,\theta)}_{\mu}} \leq C 	\left\lVert f \right\rVert_{\mathsf{W}^{l}\mathrm{L}^{(p,\theta)}_{\mu}}. 
			\end{align*}
			
			\item[(ii)] There exist $\alpha \in \mathsf{W}^{l+1}\mathrm{L}^{(p,\theta)}_{\mu}\left( \Omega; \varLambda^{k-1}\right)$ and $\beta \in \mathsf{W}^{l+1}\mathrm{L}^{(p,\theta)}_{\mu}\left( \Omega; \varLambda^{k+1}\right)$ and $h \in \mathcal{H}_{N}\left(\Omega; \varLambda^{k}\right)$ such that 
			\begin{gather*}
				f = d\alpha  + d^{\ast}\left( A\left(x\right)\beta\right) + h \qquad\text{ in } \Omega, \\
				d^{\ast}\alpha = d\beta = 0 \quad \text{ in } \Omega, \qquad 
				\nu\lrcorner \alpha  = \nu \lrcorner \left( A\left(x\right)\beta \right) = 0 \quad \text{ on } \partial\Omega.
			\end{gather*} 
			Moreover, we have the estimate 
			\begin{align*}
				\left\lVert \alpha \right\rVert_{\mathsf{W}^{l+1}\mathrm{L}^{(p,\theta)}_{\mu}} + 	\left\lVert \beta \right\rVert_{\mathsf{W}^{l+2}\mathrm{L}^{(p,\theta)}_{\mu}} + 	\left\lVert h \right\rVert_{\mathsf{W}^{l}\mathrm{L}^{(p,\theta)}_{\mu}} \leq C 	\left\lVert f \right\rVert_{\mathsf{W}^{l}\mathrm{L}^{(p,\theta)}_{\mu}}. 
			\end{align*}
		\end{itemize}
	\end{theorem}
	\begin{remark}
		By Hodge duality, each of the above cases imply their Hodge dual versions as well. 
	\end{remark}
	\begin{proof}
		We only show $(ii)$ for $l=0.$ Pick any exponent $q>1$ such that $\mathcal{L}^{\left(p, \theta\right)} \subset L^{q}.$ Using the standard Hodge decomposition in $L^{q},$ we can write 
		\begin{align*}
			f &= g + h &&\text{ in } \Omega, 
		\end{align*}
		where $h \in \mathcal{H}_{N}\left(\Omega; \varLambda^{k}\right)$ and $g \in L^{q}\left(\Omega; \varLambda^{k}\right)$ satisfies 
		\begin{align*}
			\int_{\Omega} \left\langle g, \psi \right\rangle = 0 \qquad \text{ for all } \psi \in \mathcal{H}_{N}\left(\Omega; \varLambda^{k}\right),
		\end{align*} and the estimate 
		\begin{align*}
			\left\lVert g \right\rVert_{L^{q}} + \left\lVert h \right\rVert_{L^{q}} \leq C \left\lVert f \right\rVert_{L^{q}}.
		\end{align*}
		Moreover, since harmonic fields are smooth and $L^{q}$ norm of any derivatives can be controlled by the $L^{q}$ norm of a harmonic field, we have the estimates 
		\begin{align*}
			\left\lVert h \right\rVert_{\mathrm{L}^{(p,\theta)}_{\mu}} \leq C \left\lVert h \right\rVert_{L^{q}} \leq C \left\lVert f \right\rVert_{L^{q}} \leq C 	\left\lVert f \right\rVert_{\mathrm{L}^{(p,\theta)}_{\mu}}.
		\end{align*}
		This obviously implies 
		\begin{align*}
			\left\lVert g \right\rVert_{\mathrm{L}^{(p,\theta)}_{\mu}} \leq C \left(	\left\lVert f \right\rVert_{\mathrm{L}^{(p,\theta)}_{\mu}} + 	\left\lVert h \right\rVert_{\mathrm{L}^{(p,\theta)}_{\mu}}  \right) \leq C 	\left\lVert f \right\rVert_{\mathrm{L}^{(p,\theta)}_{\mu}}.
		\end{align*}
		Now, we use Theorem \ref{generalHodgesystemtheoremnormal} and Remark \ref{remarkhodgegeneralnormal} to find $\omega \in \mathsf{W}^{2}\mathrm{L}^{(p,\theta)}_{\mu}\left( \Omega; \varLambda^{k}\right)$ that uniquely solves the system 
		\begin{align*}
			\left\lbrace \begin{gathered}
				d^{\ast} ( A (x) d\omega )  + d d^{\ast} \omega    = g  \text{ in } \Omega, \\
				\nu\lrcorner \omega  = 0 \text{  on } \partial\Omega, \\
				\nu\lrcorner  \left( A(x)d\omega \right) = 0 \text{ on } \partial\Omega. 
			\end{gathered} 
			\right.
		\end{align*}
		Setting $\alpha = d^{\ast}\omega$ and $\beta = d\omega$ completes the proof. 
	\end{proof}
	\subsection{Morrey-Lorentz estimate for Maxwell systems}
	\begin{theorem}\label{Maxwell in MorreyLorentz}
		Let $A \in C^{l+1}\left(\overline{\Omega}; \operatorname{Hom}\left(\varLambda^{k+1}\right)\right)$ and $B \in C^{l+2}\left(\overline{\Omega}; \operatorname{Hom}\left(\varLambda^{k}\right)\right)$ satisfy the Legendre condition. 
		Let  $\omega_{0} \in \mathsf{W}^{l+2}\mathrm{L}^{(p,\theta)}_{\mu}\left( \Omega; \varLambda^{k}\right),$
		$f \in \mathsf{W}^{l}\mathrm{L}^{(p,\theta)}_{\mu}\left( \Omega; \varLambda^{k}\right), $ and $g  \in \mathsf{W}^{l+1}\mathrm{L}^{(p,\theta)}_{\mu}\left( \Omega; \varLambda^{k}\right),$ and $\lambda \geq 0.$ Suppose $f,$ $g$ and $\lambda$  satisfy 
		\begin{align}
			d^{\ast} f + \lambda g = 0 \text{ and } d^{\ast} g = 0 \text{ in } \Omega . \tag{C}
		\end{align}
		\begin{itemize}
			\item[(i)] Suppose $ g \in \left(\mathcal{H}_T(\Omega;\varLambda^{k-1})\right)^{\perp}$ and if $\lambda = 0,$ assume in addition that 
			$f \in \left(\mathcal{H}_T(\Omega;\varLambda^{k})\right)^{\perp}.$  Then   
			the following boundary value problem, 
			\begin{equation} \label{problemMaxwellgeneraltan}
				\left\lbrace \begin{gathered}
					d^{\ast} ( A (x) d\omega )   = \lambda B(x)\omega + f  \text{ in } \Omega, \\
					d^{\ast} \left( B(x) \omega \right) = g \text{ in } \Omega, \\
					\nu\wedge \omega = \nu\wedge \omega_{0} \text{  on } \partial\Omega,
				\end{gathered} 
				\right. \tag{$PM_{T}$}
			\end{equation}
			has a unique solution $\omega \in \mathsf{W}^{l+2}\mathrm{L}^{(p,\theta)}_{\mu}\left( \Omega; \varLambda^{k}\right),$
			satisfying the estimates 
			\begin{align*}
				\left\lVert \omega \right\rVert_{\mathsf{W}^{l+2}\mathrm{L}^{(p,\theta)}_{\mu}} \leq c \left(  \left\lVert  \omega \right\rVert_{\mathrm{L}^{(p,\theta)}_{\mu}} + \left\lVert f\right\rVert_{\mathsf{W}^{l}\mathrm{L}^{(p,\theta)}_{\mu}}  
				+ \left\lVert g\right\rVert_{\mathsf{W}^{l}\mathrm{L+1}^{(p,\theta)}_{\mu}} + \left\lVert \omega_{0}\right\rVert_{\mathsf{W}^{l+2}\mathrm{L}^{(p,\theta)}_{\mu}}\right), 
			\end{align*}
			
			\item[(ii)]  Suppose \begin{align*}
				\nu\lrcorner g = \nu \lrcorner d^{\ast} \left(  B(x)\omega_{0} \right) \quad\text{ and }\quad \nu\lrcorner f 
				= \nu \lrcorner \left[ d^{\ast} \left( A(x)d\omega_{0} \right) - \lambda 
				B(x)\omega_{0} \right] \quad\text{ on } \partial\Omega
			\end{align*}
			and 
			$$ \int_{\Omega} \left\langle g ; \psi \right\rangle  - \int_{\partial\Omega} \left\langle \nu\lrcorner \left( B(x)\omega_{0}\right); \psi \right\rangle = 0 \qquad \text{ for all } \psi 
			\in \mathcal{H}_N(\Omega;\varLambda^{k-1}).$$
			If $\lambda =0,$ assume in addition that  
			$$\int_{\Omega} \left\langle f ; \phi \right\rangle  - \int_{\partial\Omega} \left\langle \nu\lrcorner \left( A(x)d\omega_{0} \right); \phi \right\rangle = 0 \qquad \text{ for all } \phi 
			\in \mathcal{H}_N(\Omega;\varLambda^{k}). $$ Then the following boundary value problem, 
			\begin{equation} \label{problemMaxwellgeneralnormal}
				\left\lbrace \begin{gathered}
					d^{\ast} ( A (x) d\omega )   = \lambda B(x)\omega + f  \text{ in } \Omega, \\
					d^{\ast} \left( B(x) \omega \right) = g \text{ in } \Omega, \\
					\nu\lrcorner \left( B(x) \omega \right) = \nu\lrcorner \left( B(x) \omega_{0} \right)  \text{  on } \partial\Omega, \\
					\nu\lrcorner \left( A(x) d\omega \right) = \nu\lrcorner \left( A(x) d\omega_{0} \right)  \text{  on } \partial\Omega. 
				\end{gathered} 
				\right. \tag{$PM_{N}$}
			\end{equation}
			has a unique solution $\omega \in \mathsf{W}^{l+2}\mathrm{L}^{(p,\theta)}_{\mu}\left( \Omega; \varLambda^{k}\right),$  
			satisfying the estimates 
			\begin{align*}
				\left\lVert \omega \right\rVert_{\mathsf{W}^{l+2}\mathrm{L}^{(p,\theta)}_{\mu}} \leq c \left(  \left\lVert  \omega \right\rVert_{\mathrm{L}^{(p,\theta)}_{\mu}} + \left\lVert f\right\rVert_{\mathsf{W}^{l}\mathrm{L}^{(p,\theta)}_{\mu}}  
				+ \left\lVert g\right\rVert_{\mathsf{W}^{l+1}\mathrm{L}^{(p,\theta)}_{\mu}} + \left\lVert \omega_{0}\right\rVert_{\mathsf{W}^{l+2}\mathrm{L}^{(p,\theta)}_{\mu}}\right). 
			\end{align*}
		\end{itemize}
	\end{theorem}
	\begin{remark}
		Once again, by Hodge duality, each of the above cases imply their Hodge dual versions as well. 
	\end{remark}
	\begin{proof}
		We need to show only part $(i)$, the other case is similar.
		At first, using Theorem \ref{generalHodgesystemtheorem}, we find $\alpha\in\mathsf{W}^{l+3}\mathrm{L}^{(p,\theta)}_{\mu}\left( \Omega; \varLambda^{k}\right)$ such that
		\begin{equation*}
			\left\lbrace \begin{aligned}
				d^{\ast} ( B (x) d\alpha )  + d d^{\ast}\alpha   &=  g - d^{\ast}\left(B(x) \omega_0\right)  &&\text{ in } \Omega, \\
				\nu\wedge \alpha &= 0 &&\text{  on } \partial\Omega, \\
				\nu\wedge d^{\ast} \alpha &= 0  &&\text{ on } \partial\Omega. 
			\end{aligned} 
			\right. 
		\end{equation*}
		Now it is easy to see that $\beta := d^{\ast}\alpha$ solves 
		\begin{equation*}
			\left\lbrace \begin{aligned}
				\left(d^{\ast}d + dd^{\ast}\right) \beta  &= 0  &&\text{ in } \Omega, \\
				\nu\wedge \beta &= 0 &&\text{  on } \partial\Omega, \\
				\nu\wedge d^{\ast} \beta &= 0  &&\text{ on } \partial\Omega. 
			\end{aligned} 
			\right. 
		\end{equation*}
		Now, by uniqueness of solutions to the above system, we get $\beta = d^{\ast} \theta \equiv  0$. Now, once again by Theorem \ref{generalHodgesystemtheorem}, we find $u \in \mathsf{W}^{l+2}\mathrm{L}^{(p,\theta)}_{\mu}\left( \Omega; \varLambda^{k}\right)$  solving 
		\begin{equation}\label{uequation}
			\left\lbrace \begin{aligned}
				d^{\ast} ( A (x) du )  + \left( B(x) \right)^{T}d d^{\ast}\left( B(x) u \right)   &=  \lambda B(x)u + \widetilde{f}  &&\text{ in } \Omega, \\
				\nu\wedge u &= 0 &&\text{  on } \partial\Omega. \\
				\nu\wedge d^{\ast} \left( B(x)u\right) &= 0  &&\text{ on } \partial\Omega. 
			\end{aligned} 
			\right. 
		\end{equation}
		where 
		$\widetilde{f} = f + \lambda B(x)\omega_{0} + \lambda B(x)G - d^{\ast} ( A (x) d\omega_{0} ).$ Note that if $\lambda = 0,$ 
		then $\widetilde{f} \in \left(\mathcal{H}_T(\Omega;\varLambda^{k})\right)^{\perp}$ and 
		$d^{\ast}\widetilde{f} = 0$ and thus \eqref{uequation} can always be solved for any $\lambda \geq 0.$ (see remark \ref{remarkhodgegeneraltan}).
		But this implies $v= d^{\ast} \left( B(x)u \right)$ solves the system 
		\begin{align*}
			\left\lbrace \begin{aligned}
				d^{\ast}\left( B^{T}(x) dv  \right) &= \lambda v &&\text{ in } \Omega, \\
				d^{\ast}v &=0 &&\text{ in } \Omega, \\
				\nu\wedge v &= 0  &&\text{ on } \partial\Omega. \\
			\end{aligned}\right. 
		\end{align*}
		But this implies 
		\begin{align*}
			\gamma \int_{\Omega}\left\lvert dv \right\rvert^{2} \leq \int_{\Omega}\left\langle \left( B(x)\right)^{\intercal} dv; dv \right\rangle = \lambda \int_{\Omega}\left\lvert v \right\rvert^{2}. 
		\end{align*}
		This implies $\lambda >0$ is impossible for nontrivial $v$ and if $\lambda = 0,$ $v$ must be a harmonic field. But no nontrivial harmonic field can be coexact. Hence  
		in either case, we deduce  $ v = d^{\ast} \left( B(x)u \right) \equiv 0$ in $\Omega.$ Now it is easy to check that  $\omega = \omega_{0} + u + d\alpha$ solves \eqref{problemMaxwellgeneraltan}.  
	\end{proof}
	\subsection{Morrey-Lorentz estimate for $\operatorname{div}$-$\operatorname{curl}$ systems}
	\begin{theorem}\label{divcurl system MorreyLorentz}
		Let  $A,B \in C^{l+1}(\overline{\Omega} ; \operatorname{Hom}(\varLambda^{k})),$ satisfy the 
		Legendre condition. Let $\omega_{0} \in \mathsf{W}^{l+1}\mathrm{L}^{(p,\theta)}_{\mu}\left( \Omega; \varLambda^{k}\right)$, 
		$f \in \mathsf{W}^{l}\mathrm{L}^{(p,\theta)}_{\mu}\left( \Omega; \varLambda^{k+1}\right),$ and $g \in \mathsf{W}^{l}\mathrm{L}^{(p,\theta)}_{\mu}\left( \Omega; \varLambda^{k-1}\right).$
		Then the following hold true.\smallskip
		
		\noindent \textbf{(i)} Suppose $f$ and $g$  satisfy $df = 0$, $d^{\ast} g = 0$ in $\Omega$ and 
		$ \nu\wedge d\omega_{0} = \nu\wedge f $ on $\partial\Omega,$
		and for every $\chi \in \mathcal{H}_T(\Omega;\varLambda^{k+1})$ and $\psi \in \mathcal{H}_T(\Omega;\varLambda^{k-1})$,
		\begin{equation*}
			\int_{\Omega} \langle f ; \chi \rangle - \int_{\partial\Omega} \langle \nu \wedge \omega_0 ; \chi \rangle = 0 
			\text{ and } \int_{\Omega} \langle g ; \psi \rangle = 0. 
		\end{equation*}
		Then there exists a solution $\omega \in \mathsf{W}^{l+2}\mathrm{L}^{(p,\theta)}_{\mu}\left( \Omega; \varLambda^{k}\right),$ to the following boundary value problem, 
		\begin{equation} \label{problemddeltalinear}
			\left\lbrace \begin{aligned}
				d(A(x)\omega) = f  \quad &\text{and} \quad  d^{\ast} (B(x) \omega) = g &&\text{ in } \Omega, \\
				\nu\wedge A(x)\omega &= \nu\wedge\omega_0 &&\text{  on } \partial\Omega,
			\end{aligned} 
			\right. \tag{$\mathcal{P}_{T}$}
		\end{equation}
		satisfying the estimates 
		\begin{align*}
			\left\lVert \omega \right\rVert_{\mathsf{W}^{l+1}\mathrm{L}^{(p,\theta)}_{\mu}} \leq c \left(  \left\lVert  \omega \right\rVert_{\mathrm{L}^{(p,\theta)}_{\mu}} + \left\lVert f\right\rVert_{\mathsf{W}^{l}\mathrm{L}^{(p,\theta)}_{\mu}}  + \left\lVert g\right\rVert_{\mathsf{W}^{l}\mathrm{L}^{(p,\theta)}_{\mu}} + \left\lVert \omega_{0}\right\rVert_{\mathsf{W}^{l+1}\mathrm{L}^{(p,\theta)}_{\mu}}\right), 
		\end{align*} \smallskip
		
		\noindent\textbf{(ii)} Suppose $f$ and $g$  satisfy $df = 0$, $d^{\ast} g = 0$ in $\Omega$ and $ \nu\lrcorner g  = \nu\lrcorner d^{\ast}\omega_{0}$ on $\partial\Omega,  $
		and for every $\chi \in \mathcal{H}_N(\Omega;\varLambda^{k-1})$ and $\psi \in \mathcal{H}_N(\Omega;\varLambda^{k+1})$,
		\begin{equation*}
			\int_{\Omega} \langle g ; \chi \rangle - \int_{\partial\Omega} \langle \nu \lrcorner \omega_0 ; \chi \rangle = 0 
			\text{ and } \int_{\Omega} \langle f ; \psi \rangle = 0. 
		\end{equation*}
		Then there exists a solution $\omega \in \mathsf{W}^{l+1}\mathrm{L}^{(p,\theta)}_{\mu}\left( \Omega; \varLambda^{k}\right),$ to the following boundary 
		value problem,
		\begin{equation} \label{problemddeltalinearnormal}
			\left\lbrace \begin{aligned}
				d(A(x)\omega) = f  \quad &\text{and} \quad  d^{\ast} (B(x) \omega) = g &&\text{ in } \Omega, \\
				\nu\lrcorner B(x)\omega &= \nu\lrcorner\omega_0 &&\text{  on } \partial\Omega,
			\end{aligned} 
			\right. \tag{$\mathcal{P}_{N}$}
		\end{equation}
		satisfying the estimates 
		\begin{align*}
			\left\lVert \omega \right\rVert_{\mathsf{W}^{l+1}\mathrm{L}^{(p,\theta)}_{\mu}} \leq c \left( \left\lVert  \omega \right\rVert_{\mathrm{L}^{(p,\theta)}_{\mu}} +  \left\lVert f\right\rVert_{\mathsf{W}^{l}\mathrm{L}^{(p,\theta)}_{\mu}}  + \left\lVert g\right\rVert_{\mathsf{W}^{l}\mathrm{L}^{(p,\theta)}_{\mu}} + \left\lVert \omega_{0}\right\rVert_{\mathsf{W}^{l+1}\mathrm{L}^{(p,\theta)}_{\mu}}\right), 
		\end{align*}
	\end{theorem}
	\begin{proof}
		We prove only part (ii). We use Theorem \ref{Hodge decomposition} to write 	$g - d^{\ast}\omega_{0} = d\alpha + d^{\ast}\beta + h , $	where \begin{align*}
			d^{\ast}\alpha = d\beta = 0 \text{ in } \Omega, \quad 
			\nu\lrcorner \alpha  = \nu \lrcorner \beta  = 0 \text{ on } \partial\Omega. 
		\end{align*}
		Using the hypotheses on $g$, it is easy to see that $\alpha$ and $h$ must vanish identically. Indeed, $\alpha$ satisfies 
		\begin{align*}
			\left\lbrace	\begin{aligned}
				\left( d^{\ast}d + dd^{\ast}\right) \alpha &=0 &&\text{ in } \Omega, \\
				\nu\lrcorner \alpha &= 0 &&\text{ on } \partial\Omega, \\
				\nu\lrcorner d^{\ast}\alpha &=0 &&\text{ on } \partial\Omega. 
			\end{aligned}		\right. 
		\end{align*}
		To see $h$ must vanish, we note that 
		\begin{align*}
			0 = \int_{\Omega} \langle g ; h \rangle - \int_{\partial\Omega} \langle \nu \lrcorner \omega_0 ; h \rangle = \int_{\Omega} \langle g - d^{\ast}\omega_{0} ; h \rangle = \int_{\Omega} \left\lvert h \right\rvert^{2}. 
		\end{align*}
		Now we define the matrix field  $D:=AB^{-1},$ which is clearly uniformly elliptic as well and find $\psi \in \mathsf{W}^{l+2}\mathrm{L}^{(p,\theta)}_{\mu}\left( \Omega; \varLambda^{k-1}\right)$ such that 
		\begin{align*}
			\left\lbrace \begin{aligned}
				d\left( D\left(x\right)d^{\ast}\psi \right) &= f - d \left[ D\left(x\right)\left( \beta + \omega_{0}\right)\right] &&\text{ in } \Omega, \\
				d\psi &= 0 &&\text{ in } \Omega, \\
				\nu \lrcorner \psi &= 0 &&\text{ on } \partial\Omega. 
			\end{aligned}\right. 
		\end{align*}
		Note that we can solve this system as this is the Hodge dual to  \eqref{problemMaxwellgeneraltan}. 
		Now setting $ \omega = B^{-1} ( \beta + \omega_{0} + d^{\ast}\psi )$ completes the proof. 
	\end{proof}

	\subsection{Gaffney inequality in Morrey-Lorentz spaces}
	\noindent As an immediate consequence of Theorem \ref{divcurl system MorreyLorentz}, we get the following Gaffney type inequalities in Morrey-Lorentz spaces.
	\begin{theorem}\label{gaffney} (Gaffney type inequality)
		Let $1\leq k\leq n-1$, $l\geq 0$ be integers and $1<p<\infty$, $1\leq\theta<\infty$, $0\leq\mu<n$ be real numbers. Let $\Omega\subset\mathbb{R}^n$ be open, bounded and $C^2$. Let $A\in C^{1}\left(\overline\Omega; \operatorname{Hom} (\varLambda^{k+1})\right),$ $B\in C^{1}\left(\overline\Omega; \operatorname{Hom} (\varLambda^k)\right)$ satisfy the Legendre condition. Let $u \in \mathrm{L}^{(p,\theta)}_{\mu}\left( \Omega; \varLambda^{k}\right)$,  satisfy 
		\begin{align*}
			d\left( A u\right) \in \mathrm{L}^{(p,\theta)}_{\mu}\left( \Omega; \varLambda^{k+1}\right) \qquad \text{ and } \qquad d^{\ast}\left(Bu\right) \in \mathrm{L}^{(p,\theta)}_{\mu}\left( \Omega; \varLambda^{k-1}\right).
		\end{align*} Suppose either $\nu\wedge \left( A\left(x\right)u\right)  =0$ on $\partial\Omega$ or $\nu\lrcorner\left(B(x)u\right) =0$ on $\partial\Omega$. Then $u \in \mathsf{W}^{1}\mathrm{L}^{(p,\theta)}_{\mu}\left( \Omega; \varLambda^{k}\right)$ and there exists a constant $C_p = C(\gamma, \Omega, A, B, p, \theta, \mu) >0$, such that 
		\begin{align*}
			\left\lVert u \right\rVert_{\mathsf{W}^{1}\mathrm{L}^{(p,\theta)}_{\mu}} \leq C_p \left(  \left\lVert  	d\left( Au\right)  \right\rVert_{\mathrm{L}^{(p,\theta)}_{\mu}} + \left\lVert d^{\ast}\left(Bu\right)\right\rVert_{\mathrm{L}^{(p,\theta)}_{\mu}}  
			+ \left\lVert u\right\rVert_{\mathrm{L}^{(p,\theta)}_{\mu}}\right). 
		\end{align*}
	\end{theorem}


\begin{thebibliography}{10}
	
	\bibitem{Adams_RieszPotentials}
	{\sc Adams, D.~R.}
	\newblock A note on riesz potentials.
	\newblock {\em Duke Math. J. 42}, 4 (12 1975), 765--778.
	
	\bibitem{Bennet_Sharpley_InterpolationOperators}
	{\sc Bennett, C., and Sharpley, R.}
	\newblock {\em Interpolation of operators}, vol.~129 of {\em Pure and Applied
		Mathematics}.
	\newblock Academic Press, Inc., Boston, MA, 1988.
	
	\bibitem{Burenkov_SobolevSpacesDomains}
	{\sc Burenkov, V.~I.}
	\newblock {\em Sobolev spaces on domains}, vol.~137 of {\em Teubner-Texte zur
		Mathematik [Teubner Texts in Mathematics]}.
	\newblock B. G. Teubner Verlagsgesellschaft mbH, Stuttgart, 1998.
	
	\bibitem{Calderon_SobolevSpaces}
	{\sc Calder\'{o}n, A.-P.}
	\newblock Lebesgue spaces of differentiable functions and distributions.
	\newblock In {\em Proc. {S}ympos. {P}ure {M}ath., {V}ol. {IV}}. Amer. Math.
	Soc., Providence, RI, 1961, pp.~33--49.
	
	\bibitem{Costea_SobolevLorentzCapacity}
	{\sc Costea, c.}
	\newblock Sobolev-{L}orentz capacity and its regularity in the {E}uclidean
	setting.
	\newblock {\em Ann. Acad. Sci. Fenn. Math. 44}, 1 (2019), 537--568.
	
	\bibitem{CsatoDacKneuss}
	{\sc Csat{\'o}, G., Dacorogna, B., and Kneuss, O.}
	\newblock {\em {The pullback equation for differential forms}}.
	\newblock {Progress in Nonlinear Differential Equations and their Applications,
		83}. Birkh{\"a}user/Springer, New York, 2012.
	
	\bibitem{Lamberti_et_al_BurenkovExtensionMorreySobolev}
	{\sc Fanciullo, M.~S., and Lamberti, P.~D.}
	\newblock On {B}urenkov's extension operator preserving {S}obolev-{M}orrey
	spaces on {L}ipschitz domains.
	\newblock {\em Math. Nachr. 290}, 1 (2017), 37--49.
	
	\bibitem{giaquinta-martinazzi-regularity}
	{\sc Giaquinta, M., and Martinazzi, L.}
	\newblock {\em {An introduction to the regularity theory for elliptic systems,
			harmonic maps and minimal graphs}}, second~ed., vol.~11 of {\em {Appunti.
			Scuola Normale Superiore di Pisa (Nuova Serie) [Lecture Notes. Scuola Normale
			Superiore di Pisa (New Series)]}}.
	\newblock Edizioni della Normale, Pisa, 2012.
	
	\bibitem{Hatano_MorreyLorentz}
	{\sc Hatano, N.}
	\newblock Fractional operators on {M}orrey-{L}orentz spaces and the {O}lsen
	inequality.
	\newblock {\em Math. Notes 107}, 1-2 (2020), 63--79.
	
	\bibitem{Hestenes_ExtensionbyReflection}
	{\sc Hestenes, M.~R.}
	\newblock Extension of the range of a differentiable function.
	\newblock {\em Duke Math. J. 8\/} (1941), 183--192.
	
	\bibitem{Hunt_LorentzSpaces}
	{\sc Hunt, R.~A.}
	\newblock On {$L(p,\,q)$} spaces.
	\newblock {\em Enseign. Math. (2) 12\/} (1966), 249--276.
	
	\bibitem{Koskela_et_al_MorreySObolevExtension}
	{\sc Koskela, P., Zhang, Y. R.-Y., and Zhou, Y.}
	\newblock Morrey-{S}obolev extension domains.
	\newblock {\em J. Geom. Anal. 27}, 2 (2017), 1413--1434.
	
	\bibitem{LambertiViolo_extensionMorreySobolev}
	{\sc Lamberti, P.~D., and Violo, I.~Y.}
	\newblock On stein's extension operator preserving sobolev–morrey spaces.
	\newblock {\em Mathematische Nachrichten 0}, 0 (2019), 1--15.
	
	\bibitem{Lieberman_morrey_from_Lp}
	{\sc Lieberman, G.~M.}
	\newblock A mostly elementary proof of {M}orrey space estimates for elliptic
	and parabolic equations with {VMO} coefficients.
	\newblock {\em J. Funct. Anal. 201}, 2 (2003), 457--479.
	
	\bibitem{Lorentz_LorentzSpacesIntroFirst}
	{\sc Lorentz, G.~G.}
	\newblock Some new functional spaces.
	\newblock {\em Ann. of Math. (2) 51\/} (1950), 37--55.
	
	\bibitem{Morrey_MorreySpacesintro}
	{\sc Morrey, Jr., C.~B.}
	\newblock On the solutions of quasi-linear elliptic partial differential
	equations.
	\newblock {\em Trans. Amer. Math. Soc. 43}, 1 (1938), 126--166.
	
	\bibitem{MorreyHarmonic2}
	{\sc Morrey, J. C.~B.}
	\newblock {A variational method in the theory of harmonic integrals. {II}}.
	\newblock {\em Amer. J. Math. 78\/} (1956), 137--170.
	
	\bibitem{Morrey1966}
	{\sc Morrey, J. C.~B.}
	\newblock {\em {Multiple integrals in the calculus of variations}}.
	\newblock {Die Grundlehren der mathematischen Wissenschaften, Band 130}.
	Springer-Verlag New York, Inc., New York, 1966.
	
	\bibitem{ONeil_ConvolutionLorentzSpaces}
	{\sc O'Neil, R.}
	\newblock Convolution operators and {$L(p,\,q)$} spaces.
	\newblock {\em Duke Math. J. 30\/} (1963), 129--142.
	
	\bibitem{Peetre_Lorentzembedding}
	{\sc Peetre, J.}
	\newblock Espaces d'interpolation et th\'{e}or\`eme de {S}oboleff.
	\newblock {\em Ann. Inst. Fourier (Grenoble) 16}, fasc., fasc. 1 (1966),
	279--317.
	
	\bibitem{SchwarzHodge}
	{\sc Schwarz, G.}
	\newblock {\em {Hodge decomposition---a method for solving boundary value
			problems}}, vol.~1607 of {\em {Lecture Notes in Mathematics}}.
	\newblock Springer-Verlag, Berlin, 1995.
	
	\bibitem{silthesis}
	{\sc Sil, S.}
	\newblock {Calculus of {V}ariations for {D}ifferential {F}orms, PhD Thesis}.
	\newblock {\em EPFL}, Thesis No. 7060 (2016).
	
	\bibitem{Sil_linearregularity}
	{\sc Sil, S.}
	\newblock Regularity for elliptic systems of differential forms and
	applications.
	\newblock {\em Calc. Var. Partial Differential Equations 56}, 6 (2017), 56:172.
	
	\bibitem{Tartar_Lorentzembeddings}
	{\sc Tartar, L.}
	\newblock Imbedding theorems of {S}obolev spaces into {L}orentz spaces.
	\newblock {\em Boll. Unione Mat. Ital. Sez. B Artic. Ric. Mat. (8) 1}, 3
	(1998), 479--500.
	
\end{thebibliography}
\end{document}